\documentclass[10pt,a4 paper]{article}
\usepackage{amssymb}
\usepackage{amsmath}
\usepackage{amsthm}
\usepackage{amsfonts}
\usepackage{amstext}
\usepackage{cite}
\usepackage{color}

\newtheorem{theorem}{Theorem}
\newtheorem{lemma}{Lemma}
\newtheorem{proposition}{Proposition}

\newcommand{\ds}{\displaystyle}

\newcommand{\Bo}{{\hfill $\Box$\hspace{1in}\medskip}}

\begin{document}
\title{ Continuous data assimilation for the three-dimensional Navier-Stokes-$\alpha$ model\footnote{Email: deboraalbanez@utfpr.edu.br (D.A.F. Albanez), hlopes@im.ufrj.br (H.J. Nussenzveig Lopes) titi@math.tamu.edu (E.S. Titi), }}
\author{  D\'ebora A.F. Albanez$^{\lowercase{a}}$,  Helena J. Nussenzveig Lopes$^{\lowercase{b}}$, \\ Edriss S. Titi$^{c,d}$, }
\date{August 21, 2014}
\maketitle

\begin{center}

$^a$ Departamento de Matem{\'a}tica \; -- \;
Universidade Tecnol{\'o}gica Federal do Paran{\'a}, \\
86300-000 Corn{\'e}lio Proc{\'o}pio, PR -- Brasil.\\
$^b$ Instituto de Matem\'atica \; -- \;
Universidade Federal do Rio de Janeiro, \\
Cidade Universit\'aria -- Ilha do Fund\~ao, \;\;
Caixa Postal 68530,\\
21941-909 Rio de Janeiro, RJ -- Brasil.\\
$c$ Department of Computer Science and Applied Mathematics\\
Weizmann Institute of Science,
Rehovot, 76100, Israel.\\
$d$ Department of Mathematics, Texas A\&M University, 3368 TAMU \\
                 College Station, TX 77843-3368, USA.\\

\end{center}

\begin{abstract}

Motivated by the presence of a finite number of determining parameters (degrees of freedom) such as modes, nodes and local spatial averages for dissipative dynamical systems, we present a continuous data assimilation algorithm for the three-dimensional Navier-Stokes-$\alpha$ model. This algorithm consists of introducing a nudging process through general type of approximation interpolation operator (that is constructed from observational measurements) that synchronizes the large spatial scales of the approximate solutions with those of the  unknown solutions the Navier-Stokes-$\alpha$ equations that corresponds to these measurements. Our main result provides conditions on the finite-dimensional spatial resolution of the collected data, sufficient to guarantee that the approximating solution, that is obtained from these collected data, converges to the unkown reference solution (physical reality) over time. These conditions are given in terms of some physical parameters, such as kinematic viscosity, the size of the domain and the forcing term.


\end{abstract}
\numberwithin{equation}{section}

\numberwithin{equation}{section}

%
%\section{Momentum equation and equivalent form}

\section{Introduction}

The method of inserting observational measurements obtained from some physical system into the theorical model of this system, as the latter is being integrated in time, is called continuous data assimilation. It was first proposed in \cite{Daley} for atmospheric predictions, such as weather forecasting.
In general, producing accurate information of the true state of the atmosphere or a fluid in a given time is not possible, because the observational measurements are given with low spatial and/or temporal resolution.

 In order to obtain a good asymptotic approximation of the real physical state $u(t)$, the classical method of continuous data assimilation requires to separate the fast and slow parts of the solution, before inserting the measured data into the model. An application of this method of separation for 2D Navier-Stokes equations is given in \cite{OlsonTiti2003}, where the authors computed an approximation of the fast part of the solution, using an initial guess of the high modes of the exact solution. Later,  this same technique was implemented in \cite{Korn} to the 3D Navier-Stokes-$\alpha$ model.

The method we propose here for the three-dimensional Navier-Stokes-$\alpha$ model was firstly developed  in \cite{AzouaniOlsonTiti} for the two-dimensional Navier-Stokes equations, where the observational measurements are  directly inserted into the model in a way that overcomes  the difficulties coming from the fact that the discrete observations may not be  elements of a suitable space of solutions, for instance. Notably, the algorithm developed in \cite{AzouaniOlsonTiti} was inspired by the feedback control algorithm introduced  in \cite{AzouaniTiti} for stabilizing unstable solutions of dissipative partial differential equations by employing  only finitely many observables and controllers.

The aim of this work is to establish  sufficient  conditions on the spatial resolution of the observational data, and on the  relaxation (nudging) parameter $\mu>0$, that  will  guarantee the existence of an approximate solution of the real state over time. The advantage of this method is that our initial data can be chosen to be arbitrary.

Suppose that the evolution of $u$ is governed by the three-dimensional Navier-Stokes-$\alpha$ equations (cf. \cite{TitiHolmFoias}), subject to periodic boundary conditions, on  $\Omega=[0,L]^{3}$:
\begin{equation}\label{evolutionU}
\left\{\begin{array}{l}
\ds\frac{\partial v}{\partial t}-\nu\Delta v
-u\times(\nabla\times v)+\nabla p=f,  \\
\mbox{div}\,\  u=0,
\end{array}\right.
\end{equation}
on the interval $[0,T]$, where the initial data $u_{0}$ is unkown. Here $u=u(x,t)$ represents the velocity of the fluid, called the {\it filtered velocity} and $v=u-\alpha^{2}\Delta u$; $\alpha>0$ is a scale parameter with dimension of length.

 Let $I_{h}(u(t))$ represent the observational measurements at a spatial resolution of size $h$, for $t\in[0,T]$. The algorithm we use to construct an approximation $w(t)$ of $u(t)$ from the measured data consists in solving
$$\left\{\begin{array}{lcl}\label{evolutionUassimilada}
\ds\frac{\partial z}{\partial t}-\nu\Delta z
-w \!\!\!\!&\times &\!\!\!(\nabla  \times z)+\nabla p  = f \\
& -&\!\!\!\! \mu(I_{h}(w)-I_{h}(u))+\mu\alpha^{2}\Delta(I_{h}(w)-I_{h}(u)),\\
w(0)=w_{0},& & \nonumber
\end{array}\right.$$
with $z=w-\alpha^{2}\Delta w$, while  $w(0)=w_{0}$ is chosen to be arbitrary, yet in a suitable space.

This method requires that the observational measurements $I_{h}(u)$ be given as linear interpolant observables satisfying either $I_{h}:\dot{H^{1}}(\Omega)\rightarrow \dot{L^{2}}(\Omega)$ and, for some constant $c_1>0$,
\begin{equation}\label{caso1}
\|\varphi-I_{h}\varphi\|_{L^{2}(\Omega)}^{2}\leq c_{1}^{2}h^{2}\|\nabla\varphi\|^{2}, \,\,\, \mbox{for every}\,\,\, \varphi\in H^{1}(\Omega),
\end{equation}
or $I_{h}:\dot{H^{2}}(\Omega)\rightarrow \dot{L^{2}}(\Omega)$ and, for some constant $c_2>0$,
\begin{equation}\label{caso2}
\|\varphi-I_{h}\varphi\|_{L^{2}(\Omega)}^{2}\leq c_{2}^{2}h^{2}\|\nabla\varphi\|^{2}_{L^{2}(\Omega)}+{c_{2}^{2}}h^{4}\|\varphi\|^{2}_{H^{2}(\Omega)}, \,\,\, \mbox{for every}\,\,\, \varphi\in H^{1}(\Omega).
\end{equation}

One example of an interpolant operator satisfying (\ref{caso1}) and two examples of interpolant operators satisfying (\ref{caso2}) are given in Section 4. Motivated by the fact that dissipative evolution equations posses finitely many determining parameters (degrees of freedom) these  examples of interpolant operators  include determining modes (see \cite{FoiasProdi},\cite{JonesTiti2} \cite{OlsonTiti2003},\cite{OlsonTiti2008}), determining nodes (see \cite{FoiasTemam84},\cite{FoiasTiti},\cite{JonesTiti2}), determining finite volume elements (see \cite{FoiasTiti}, \cite{JonesTiti},\cite{JonesTiti2}) and finite-elements method (cf. \cite{CockburnJonesTiti}, and references therein).

Recently,  similar ideas to the data assimilation algorithms introduced in
\cite{OlsonTiti2003} and \cite{AzouaniOlsonTiti} have been implemented in \cite{FoiasJolly} \cite{FoiasJollyKravchenkoTiti} to show that the long-time dynamics of
the two-dimensional Navier--Stokes equations can be imbedded in an
infinite-dimensional dynamical system that is induced by an ordinary
differential equation in appropriate Banach space, named {\it determining form}, which is governed by a globally Lipschitz vector field. In particular, it is shown in \cite{FoiasJollyKravchenkoTiti} that solutions of the determining form converge to one of its  steady states, and that there is one to one correspondence between these steady states and the trajectories on the global attractor of the Navier--Stokes equations.
Moreover, it is worth
mentioning that the method of data assimilation studied here  can be equally applied for designing signal synchronization algorithms. Furthermore, we observe that most recently the data assimilation algorithm  that was introduced   in  \cite{AzouaniOlsonTiti} has been extended in  \cite{BessaihOlsonTiti}   to cover the case
where the observational measurements  are contaminated with stochastic random errors.

 In this paper, apply  this new data assimilation algorithm \cite{AzouaniOlsonTiti} to 3D Navier-Stokes-$\alpha$ model, in the absence of measurement errors. The  paper is organized as follows: first, we recall the functional setting of three-dimensional Navier-Stokes-$\alpha$ equations needed to develop our method of continuous data assimilation. Subsequently, we present this method and the results of well-posedness for the new data assimilation equations. Later, in section 3, we state and prove our main result, in which we give conditions under which the approximate solutions, obtained by this algorithm, converge to the solution of NS-$\alpha$ equations. Finally, in section 4, we present some examples of interpolant operators.

\section{Preliminaries and Results}

 In this section, we review some  basic facts and the functional setting of the three-dimensional Navier-Stokes-$\alpha$ equations that will be used in this paper.

Let $\Omega=[0,L]^{3}$ be a periodic box, for some $L>0$ fixed. We denote by $\mathcal{V}$ the set
 of all vector valued trigonometric polynomials defined in $\Omega$, which  are divergence-free and have average zero. Denote also by $H$ and $V$ the closure of
$\mathcal{V}$ in  the $(L^{2}(\Omega))^{3}$ and $(H^{1}(\Omega))^{3}$, respectively. The spaces $H$ and $V$ are Hilbert spaces with inner products given by
$$(u,v)=\ds\int_{\Omega}u(x)\cdot v(x)dx\,\,\, \mbox{and}\,\,\, ((u,v))=\ds\sum_{i=1}^{3}\int_{\Omega}\ds\frac{\partial u_{i}}{\partial x_{j}}\ds\frac{\partial u_{i}}{\partial x_{j}}dx, $$
respectively. Denote the norms of $H$ and $V$ by
$$|u|=(u,u)^{\frac{1}{2}}\,\,\, \mbox{and}\,\,\, \|u\|=((u,u))^{\frac{1}{2}}.$$

If $Z\subset L^{1}(\Omega)$, then we set $\dot{Z}=\{\varphi\in Z,\,\, \mbox{such that}\,\, \ds\int_{\Omega}\varphi(x)dx=0\}$.

We denote $\mathcal{P}:\dot{L^{2}}(\Omega)\rightarrow H$ the Leray projector, and by $A=-\mathcal{P}\Delta$ the Stokes operator, with domain $D(A)=H^{2}(\Omega)\cap V$. In the case of periodic boundary conditions, $A=-\Delta|_{D(A)}$. The Stokes operator is a self-adjoint positive operator, whose inverse $A^{-1}$ is a compact operator in $H$. Hence there exists a complete orthonormal set of eigenfunctions $\{\phi_{i}\}_{i=1}^{\infty}\subset H$, such that $A\phi_{i}=\lambda_{i}\phi_{i}$ and $0<\lambda_{1}\leq\lambda_{2}\leq\dots$. Let us denote $\lambda_{1}$ the first eigenvalue, i.e., $\lambda_{1}=(2\pi/L)^{2}$. We have the following versions of Poincar{\'e} inequalities: for all $u\in V$ and $v\in D(A)$,
$$|u|^{2}\leq\lambda_{1}^{-1}\|u\|^{2}\,\,\, \mbox{and}\,\,\, \|v\|^{2}\leq\lambda_{1}^{-1}|Av|^{2}. $$

Following the classical notation for Navier-Stokes equations, we write $B=B(u,v)$, $B:V\times V\rightarrow V'$, for the bilinear operator $B(u,v)=\mathcal{P}[(u\cdot\nabla)v]$.

Let us denote, for every $u,v\in \mathcal{V}$, $\widetilde{B}(u,v)=-\mathcal{P}(u\times(\nabla\times v))$. The operator $\widetilde{B}$ can be extended continuously from $V\times V$ with values in $V'$, and in particular it satisfies the following properties (see \cite{TitiHolmFoias}):
for every $u\in H,v\in V$ and $w\in D(A)$,
\begin{equation}\label{desiBtil5}
|\langle\widetilde{B}(u,v),w\rangle_{D(A)',D(A)}|\leq k_{1}|u|\, \|v\|\, \|w\|^{\frac{1}{2}}|Aw|^{\frac{1}{2}}.
\end{equation}
Also, for every $u\in V,v\in H$ and $w\in D(A)$, we have
\begin{equation}\label{desiBtil7}
|\langle\widetilde{B}(u,v),w\rangle_{D(A)',D(A)}|\leq k_{2}(|u|^{\frac{1}{2}}\|u\|^{\frac{1}{2}}|v|\,
|Aw|+|v|\, \|u\|\, \|w\|^{\frac{1}{2}}|Aw|^{\frac{1}{2}}).
\end{equation}
In addition, for every $u\in D(A),v\in H$ and $w\in V,$
\begin{equation}\label{desiBtil8}
|\langle\widetilde{B}(u,v),w\rangle_{V',V}|\leq k_{3}(\|u\|^{\frac{1}{2}}|Au|^{\frac{1}{2}}|v|\, \|w\|+|Au|\, |v|\, |w|^{\frac{1}{2}}\|w\|^{\frac{1}{2}}).
\end{equation}
The positive constants denoted by $k_{1},k_{2}$ and $k_{3}$ are scale-invariant. Furthermore, for every $u,v,w\in V$,
\begin{equation}\label{desiBtil3}
\langle\widetilde{B}(u,v),w\rangle_{V',V}=-\langle\widetilde{B}(w,v),u\rangle
_{V',V}
\end{equation}
and in particular, for every $u,v\in V$,
\begin{equation}\label{desiBtil4}
\langle\widetilde{B}(u,v),u\rangle_{V',V}=0.
\end{equation}

With the above notation, we write the incompressible three-dimensional Navier-Stokes-$\alpha$ (NS-$\alpha$) in functional form as
\begin{equation}\label{versaocompridadaBtil}
\ds\frac{dv}{dt}+\nu Av+\widetilde{B}(u,v)=f,
\end{equation}
with $v(t)=u(t)+\alpha^{2}Au(t)$, and initial condition $u(0)=u_{0}\in V$. We assume $f\in H$ is time independent so that $\mathcal{P}f=f$. Equation (\ref{versaocompridadaBtil}) is globally well-posedness, as shown in \cite{TitiHolmFoias}:
\begin{theorem}[Global existence and uniqueness for NS-$\alpha$]\label{teorema2002}
Let $f\in H$ and $u_{0}\in V$. Then for any $T>0$, system $(\ref{versaocompridadaBtil})$ has a unique regular solution that satisfies:
$$u\in C([0,T];V)\cap L^{2}([0,T];D(A))\,\,\, \mbox{and}\,\,\, \ds\frac{du}{dt}\in L^{2}([0,T];H).$$
\end{theorem}

We state now estimates of the solutions $u$ of (\ref{versaocompridadaBtil}) that will be needed later in our analysis. These estimates appear in \cite{TitiHolmFoias}, in the proof of Theorem \ref{teorema2002}.
\begin{proposition}\label{estdosatratores}
 Fix $T>0$. Let $G$ be the Grashoff number $G=\ds\frac{|f|}{\nu^{2}\lambda_{1}^{3/4}}$ and suppose that $u$ is the solution given by Theorem \ref{teorema2002}. Then there exists a time $t_{0}$, which depends on $u_{0}$, such that for $t\geq t_{0}>0$ we have
\begin{equation}\label{H1NSalpha}
\|u(t)\|^{2}\leq \frac{2G^{2}\nu^{2}}{\lambda_{1}^{1/2}\alpha^{2}}.
\end{equation}
Moreover,
\begin{equation}\label{EstimAu}
\ds\int_{t}^{t+T}(\|u(s)\|^{2}+\alpha^{2}|Au(s)|^{2})ds\leq(2+\nu\lambda_{1}T)\frac{\nu G^{2}}{\lambda_{1}^{1/2}}
\end{equation}
\end{proposition}

We present now the continuous data assimilation algorithm for the incompressible three-dimensional NS-$\alpha$ equations. Let $u$ be a regular solution of (\ref{versaocompridadaBtil}) given by Theorem \ref{teorema2002} and let $I_{h}$ be a finite rank interpolation operator satisfying either (\ref{caso1}) or (\ref{caso2}). Our aim is to recover $u$ from observational measurements $I_{h}(u(t))$, that have been measured for times $t\in[0,T]$. The approximating solutions $w$ with initial condition $w_{0}\in V$, chosen arbitrarily, will be given by the solutions of the system
\begin{equation}\label{interdentrodo3}
\left\{\begin{array}{lcl}
\ds\frac{\partial }{\partial t}(w-\alpha^{2}\Delta w)-\nu\Delta(w\!\!\!\!\!\!\! &-\!\!\!\!\!\!\! &\!\! \alpha^{2}
\Delta w)-w\times(\nabla\times (w-\alpha^{2}\Delta w))+\nabla p  \\
& = &\!\!\! f-\mu(I_{h}w-I_{h}u)+\mu\alpha^{2}\Delta(I_{h}w-I_{h}u),\\
 \mbox{div}\,\  w=0,&&
\end{array}\right.
\end{equation}
on the interval $[0,T]$. Using the Leray projector, the above system is equivalent to

\begin{equation}\label{eqcomIsoP}
 \left\{\begin{array}{ll}
\ds\frac{d}{dt}(w+\alpha^{2}Aw)+\nu A(w\,\,+&\!\alpha^{2}Aw) +
  \widetilde{B}(w,w+\alpha^{2}Aw)  \\
&\!\!=f-\mu\mathcal{P}(I_{h}w-I_{h}u)+\mu\alpha^{2}\mathcal{P}\Delta(I_{h}w-I_{h}u), \\
 \mbox{div}\,\  w=0.
\end{array}\right.
\end{equation}

Furthermore, inequalities (\ref{caso1}) and (\ref{caso2}) imply
\begin{equation}\label{become1}
|\mathcal{P}(\varphi-I_{h}\varphi)|^{2}\leq c_{1}^{2}h^{2}\|\varphi\|^{2}
\end{equation}
for all $\varphi \in V$ and
\begin{equation}\label{become2}
|\mathcal{P}(\varphi-I_{h}\varphi)|^{2}\leq c_{2}^{2}h^{2}\|\varphi\|^{2}+c_{2}^{2}h^{4}|A\varphi|^{2}.
\end{equation}
for $\varphi\in D(A)$. We will present later specific exampes of $I_{h}$ (see also \cite{AzouaniOlsonTiti}).

\begin{lemma}\label{lema1}
Suppose $\varphi\in \dot{H^{2}}(\Omega)$. Then $\mathcal{P}\varphi\in \dot{H^{2}}(\Omega)$ and $-\mathcal{P}\Delta\varphi=-\mathcal{P}\Delta\mathcal{P}\varphi=A\mathcal{P}\varphi$. 

\end{lemma}
{\it Proof.} If $\varphi\in \dot{H^{2}}(\Omega)$, by the Helmholtz decomposition (see \cite{Majda}), there exists a unique $\psi\in V$ and $p\in \dot{H^{1}}(\Omega)$ such that $\varphi=\psi+\nabla p$, with div $\psi=0$ and $\mathcal{P}\varphi=\psi$. Moreover, we also have $\Delta p=\mbox{div}\, \varphi\in\dot{H^{1}}(\Omega)$, and it follows that $p\in \dot{H^{3}}(\Omega)$. Since $\varphi\in\dot{H^{2}(\Omega)}$, we conclude that $\psi\in H^{2}(\Omega)$. On the other hand,
$$-\Delta\varphi=-\Delta\psi-\nabla(\Delta p),$$
and consequently,
\begin{equation}\label{comut}
-\mathcal{P}\Delta\varphi=-\Delta\psi=-\Delta\mathcal{P}\varphi.
\end{equation}
This also implies that
$$A\varphi=-\mathcal{P}(\Delta\varphi)=-\mathcal{P}^{2}(\Delta\varphi)=-\mathcal{P}\Delta(\mathcal{P}\varphi)=A\mathcal{P}\varphi,$$
as desired.

\Bo

Using Lemma \ref{lema1}, system (\ref{eqcomIsoP}) is equivalent to
\begin{equation}\label{eqcomInaexis}
 \left\{\begin{array}{ll}
\ds\frac{d}{dt}(w+\alpha^{2}Aw)+\nu A(w\,\,+&\!\!\!\!\alpha^{2}Aw) +
  \widetilde{B}(w,w+\alpha^{2}Aw)  \\
&\!\!=f -\mu(I+\alpha^{2}A)\mathcal{P}(I_{h}(w)-I_{h}(u)), \\
 \mbox{div}\,\  w=0,
\end{array}\right.
\end{equation}
on the interval $[0,T]$, with $w(0)=w_{0}\in V$.

Next, we show that the data assimilation equations (\ref{eqcomInaexis}) are well-posed for both cases of interpolant operators $I_{h}$: those satisfying (\ref{become1}) and those satisfying (\ref{become2}).

\begin{theorem}\label{teoexismeu}
Let $f\in H,w_{0}\in V$ and $\mu>0$ be given. Suppose that $I_{h}$ satisfies
(\ref{caso1}) (and hence $(\ref{become1})$) and that $\mu c_{1}^{2}h^{2}<\ds\frac{\nu}{2}$, where
$c_{1}>0$ is the constant given in $(\ref{become1})$. Let $u$ be the solution of NS-$\alpha$ equations with initial data $u(0)=u_{0}\in V$, ensured by Theorem \ref{teorema2002}. Then the continuous data assimilation algorithm, $(\ref{eqcomInaexis})$, has a
 regular solution $w$ that satisfies
 \begin{equation}\label{espacoNS}
 w\in C([0,T];V)\cap L^{2}([0,T];D(A))\,\,\, \mbox{and}\,\,\, \ds\frac{dw}{dt}\in L^{2}([0,T];H),
 \end{equation}
 for any $T>0$.
\end{theorem}
{\it Proof.}
First note that system (\ref{eqcomInaexis}) is equivalent to
\begin{equation}\label{inversa}
\ds\frac{d}{dt}w+\nu Aw+(I+\alpha^{2}A)^{-1}\widetilde{B}(w,z)=(I+\alpha^{2}
A)^{-1}f-\mu\mathcal{P}(I_{h}(w)-I_{h}(u)),
\end{equation}
where $z=(I+\alpha^{2}A)w$. Define
$$\overline{f}(s)=(I+\alpha^{2}A)^{-1}f+\mu\mathcal{P} I_{h}u(s).$$
Thanks to (\ref{become1}), then for all $s\in[0,T]$, we have
\begin{equation}
|\mathcal{P}I_{h}(u(s))| \leq  |\mathcal{P}(u(s)-I_{h}(u(s)))|+|u(s)| \leq c_{1}h\|u(s)\|+|u(s)|.
\end{equation}
Since the Navier-Stokes-$\alpha$ solution satisfies $u\in C([0,T];V)$, we conclude that $I_{h}(u)\in C([0,T];H)$. Moreover, we have
$$ |\overline{f}| \leq |f|+\mu c_{1}h\|u\|+\mu|u|,$$
and therefore $\overline{f}\in C([0,T];H)$, i.e., there exists  a constant $M$, that might depend on $T$, such that
$|\overline{f}|<M$, for every $t\in[0,T]$.

The purpose now is to establish the global existence of solutions to (\ref{eqcomInaexis}). For that, we use the Faedo-Galerkin method. Let $H_{m}=$ span$\{\phi_{1},\dots,\phi_{m}\}$, where $A\phi_{j}=\lambda_{j}\phi_{j}$. We denote by $P_{m}$ the orthogonal projection from $H$ onto $H_{m}$. Let $w_{m}\in H_{m}$ satisfy the finite-dimensional Faedo-Galerkin system of ordinary differential equations:
\begin{equation}\label{eqProjetadaGaler}
\left\{\begin{array}{l}
\ds\frac{dw_{m}}{dt}+\nu Aw_{m}+P_{m}(I+\alpha^{2}A)^{-1}\widetilde{B}(w_{m},z_{m})=P_{m}\overline{f}
-\mu P_{m}\mathcal{P}I_{h}(w_{m})\\
w_{m}(0)=P_{m}w_{0}.
\end{array}\right.
\end{equation}

Since system (\ref{eqProjetadaGaler}) has a quadratic non-linearity, therefore it is locally Lipschitz and as a result it has a unique short time solution. The next step is to prove that the solution is uniformly bounded in time and $m$; and thereby we shall ensure the global existence in time of $w_{m}$ for all $m$.

Denote by $[0,T_{m}^{\max})$ the maximal interval of existence for (\ref{eqProjetadaGaler}). Our goal is to show that $T_{m}^{\max}=T$. Focusing on $[0,T_{m}^{\max})$, and taking the inner product of $(\ref{eqProjetadaGaler})$ with $w_{m}+\alpha^{2}Aw_{m}$  we have

\begin{eqnarray}\label{deadd}
\ds\frac{1}{2}\frac{d}{dt}(|w_{m}|^{2}+\alpha^{2}\|w_{m}\|^{2})
 +&&\!\!\!\!\!\!\!\!\!\!\!\!\nu(\|w_{m}\|^{2}+\alpha^{2}|Aw_{m}|^{2})  \nonumber\\
&+& (P_{m}(I+\alpha^{2}A)^{-1}
\widetilde{B}(w_{m},z_{m}),w_{m}+\alpha^{2}Aw_{m}) \nonumber\\ \nonumber
&= &(P_{m}\overline{f},w_{m}+\alpha^{2}Aw_{m})-\mu
(I_{h}(w_{m}),w_{m}+\alpha^{2}Aw_{m}). \nonumber
\end{eqnarray}
Using the fact that $(I+\alpha^{2}A)$ is self-adjoint and property (\ref{desiBtil4}), we obtain

\begin{eqnarray}
\ds\frac{1}{2}\frac{d}{dt}(|w_{m}|^{2}+\alpha^{2}\|w_{m}\|^{2})
+  \nu(\|w_{m}\|^{2}+&&\!\!\!\!\!\!\!\!\!\!\!\! \alpha^{2}|Aw_{m}|^{2})
=(P_{m}\overline{f},w_{m})+\alpha^{2}(P_{m}\overline{f},Aw_{m}) \nonumber\\
&-&\!\!\!\mu(I_{h}(w_{m}),w_{m})-\mu\alpha^{2}(I_{h}(w_{m}),Aw_{m}). \nonumber  \nonumber
\end{eqnarray}
By Young's inequality,
\begin{eqnarray}
\ds\frac{1}{2}\frac{d}{dt}(|w_{m}|^{2}+\alpha^{2}\|w_{m}\|^{2})
+ &&\!\!\!\!\!\!\!\!\!\!\!\! \nu(\|w_{m}\|^{2}+\alpha^{2}|Aw_{m}|^{2})
\leq\left(\ds\frac{1}{\mu}+\ds\frac{\alpha^{2}}{\nu}\right)|\overline{f}|^{2} \nonumber \\
&+& \ds\frac{\mu}{4}|w_{m}|^{2}+\ds\frac{\nu}{4}\alpha^{2}|Aw_{m}|^{2} \label{iniciooutraexist}\\
&+&\mu(w_{m}-I_{h}(w_{m}),w_{m})-\mu|w_{m}|^{2} \nonumber \\
&+&\mu\alpha^{2}(w_{m}-I_{h}(w_{m}),Aw_{m})-\mu\alpha^{2}\|w_{m}\|^{2}.
\nonumber  \nonumber
\end{eqnarray}

Using condition $(\ref{become1})$, we have
\begin{eqnarray}
\ds\frac{1}{2}\frac{d}{dt}(|w_{m}|^{2}+\alpha^{2}\|w_{m}\|^{2})
+ &&\!\!\!\!\!\!\!\!\!\!\!\! \nu(\|w_{m}\|^{2}+\alpha^{2}|Aw_{m}|^{2})
\leq\left(\ds\frac{1}{\mu}+\ds\frac{\alpha^{2}}{\nu}\right)|\overline{f}|^{2} \nonumber\\
&+&\!\! \ds\frac{\mu}{4}|w_{m}|^{2}+\ds\frac{\nu}{4}\alpha^{2}|Aw_{m}|^{2}
 \label{umdosexis} \\
& + &\!\! \mu c_{1}h\|w_{m}\|\, |w_{m}|-\mu|w_{m}|^{2} \nonumber\\
&+ &\!\! \mu\alpha^{2}c_{1}h\|w_{m}\|\, |Aw_{m}|-
\mu\alpha^{2}\|w_{m}\|^{2}.   \nonumber
\nonumber
\end{eqnarray}

By Young's inequality again and the hypothesis that  $h$ is sufficiently small so that $\mu c_{1}^{2}h^{2}<\ds\frac{\nu}{2}$, it follows that
\begin{eqnarray}
\ds\frac{1}{2}\ds\frac{d}{dt}(|w_{m}|^{2}+
&&\!\!\!\!\!\!\!\!\!\!\!\!
\alpha^{2}\|w_{m}\|^{2})
+ \nu(\|w_{m}\|^{2}+\alpha^{2}|Aw_{m}|^{2})  \nonumber\\
&\leq &\left(\ds\frac{1}{\mu}+\ds\frac{\alpha^{2}}{\nu}\right)|\overline{f}|^{2}
+\ds\frac{\mu}{4}|w_{m}|^{2}+\frac{\nu}{4}\alpha^{2}|Aw_{m}|^{2}\nonumber\\
& + &\ds\frac{\mu}{4}|w_{m}|^{2}+\frac{\nu}{2}\|w_{m}\|^{2}-\mu|w_{m}|^{2} \nonumber\\
& + & \frac{\mu}{2}\alpha^{2}\|w_{m}\|^{2}+\frac{\nu} {4}\alpha^{2}|Aw_{m}|^{2}-\mu\alpha^{2}\|w_{m}\|^{2}.\nonumber
\nonumber
\end{eqnarray}

Therefore
\begin{eqnarray}
\ds\frac{1}{2}\ds\frac{d}{dt}(|w_{m}|^{2}+
&&\!\!\!\!\!\!\!\!\!\!\!\!
\alpha^{2}\|w_{m}\|^{2})
+ \frac{\nu}{2}(\|w_{m}\|^{2}+\alpha^{2}|Aw_{m}|^{2})
\label{paraintegrarnovo}  \\
&\leq &\left(\ds\frac{1}{\mu}+\ds\frac{\alpha^{2}}{\nu}\right)|\overline{f}|^{2}
-\ds\frac{\mu}{2}(|w_{m}|^{2}+\alpha^{2}\|w_{m}\|^{2}).\nonumber
\end{eqnarray}

 By Poincar{\'e} and Gronwall's inequality we conclude that for all $t\in[0,T_{m}^{\max})$,
\begin{equation}\label{M1primeiro}
|w_{m}(t)|^{2}+\alpha^{2}\|w_{m}(t)\|^{2}\leq(|w_{0}|^{2}+\alpha^{2}\|w_{0}\|
^{2})+2\left(\frac{1}{\mu}+\frac{\alpha^{2}}{\nu}\right)
\frac{M}{\nu\lambda_{1}+\mu}=: M_{1}.
\end{equation}

Since the right-hand side of (\ref{M1primeiro}) is bounded, then $T_{m}^{\max}=T$, otherwise we can extend the solution beyond $T_{m}^{\max}$, which contradicts the definition of $T_{m}^{\max}$.

The estimate above is uniform in $m$ and $t$, and therefore we have the global existence of $w_{m}$ in time and also
\begin{equation}\label{limLinfinito}
\|w_{m}\|^{2}_{L^{\infty}([0,T];V)}\leq\frac{M_{1}}{\alpha^{2}}\,\,\,\,\,
\mbox{and}\,\,\,\,\,\|z_{m}\|^{2}_{L^{\infty}([0,T];V')}\leq M_{1}.
\end{equation}

Additionally, from $(\ref{paraintegrarnovo})$ we get
$$\begin{array}{rl}
|w_{m}(t)|^{2}+\alpha^{2}\|w_{m}(t)\|^{2}+\nu\ds\int_{0}^{t}
(\|w_{m}(s)\|^{2}& \!\!\!\!  +  \alpha^{2}|Aw_{m}(s)|^{2})ds\\
& \leq |w_{m}(0)|^{2}+\alpha^{2}
\|w_{m}(0)\|^{2}+2\left(\ds\frac{1}{\mu}+\ds\frac{\alpha^{2}}{\nu}\right)Mt,
\end{array}$$
and it follows that
\begin{equation}\label{estimatDa}
\ds\int_{0}^{t}
(\|w_{m}(s)\|^{2}+\alpha^{2}|Aw_{m}(s)|^{2})ds\leq\frac{1}{\nu}(|w_{0}|^{2}+
\alpha^{2}\|w_{0}\|^{2})+2\left(\frac{1}{\mu}+\frac{\alpha^{2}}{\nu}\right)
\frac{MT}{\nu}=:M_{2}(T).
\end{equation}

Hence,
\begin{equation}\label{limL2DA}
\|w_{m}\|^{2}_{L^{2}([0,T];D(A))}\leq\frac{M_{2}(T)}{\alpha^{2}}\,\,\,\,\,
\mbox{and}\,\,\,\,\,\|z_{m}\|^{2}_{L^{2}([0,T];H)}\leq M_{2}(T).
\end{equation}

Note that from $(\ref{estimatDa})$ we also obtain
\begin{equation}\label{M2}
\|w_{m}\|^{2}_{L^{2}([0,T];V)}\leq M_{2}(T).
\end{equation}

Now we establish uniform estimates in $m$ for the derivatives $\ds\frac{dw_{m}(t)}{dt}$ and
 $\ds\frac{dz_{m}(t)}{dt}$. Returning to equation
\begin{eqnarray}
\ds\frac{dz_{m}(t)}{dt}+\nu Az_{m}(t)\!\!\!\!&\!\!\!\!+\!\!\!\! & \!\!\!\! P_{m}\widetilde{B}(w_{m},z_{m})=P_{m}f\label{estpraderivada}\\
& +& \!\!\!\! \mathcal{P}I_{h}(u(t))+P_{m}A\mathcal{P}I_{h}(u(s))-\mu P_{m}\mathcal{P}I_{h}(w_{m}(t))-\mu\alpha^{2}P_{m}A\mathcal{P}I_{h}(w_{m}(t)),\nonumber
\nonumber
\end{eqnarray}
and we shall estimate $\ds\frac{dz_{m}(t)}{dt}$ in $L^{2}([0,T];D(A)')$. Note that by (\ref{desiBtil7}),
\begin{eqnarray}
|\widetilde{B}(w_{m}(t),
z_{m}(t))|
  & \leq & k_{2}(|w_{m}(t)|^{1/2}\|w_{m}(t)\|^{1/2}|z_{m}(t)|+\lambda_{1}^{-1/4}
|z_{m}(t)|\|w_{m}(t)\|)\nonumber\\
 & \leq & 2k_{2}\lambda_{1}^{-1/4}|z_{m}(t)|\|w_{m}(t)\|.\nonumber
 \nonumber
\end{eqnarray}
Consequently, and thanks to (\ref{limLinfinito}) and (\ref{limL2DA}), we have
\begin{eqnarray}
\|\widetilde{B}(w_{m},z_{m})\|^{2}_{L^{2}([0,T];H)}=
&&\!\!\!\!\!\!\!\!\!\!\! \ds\int_{0}^{T}|
\widetilde{B}(w_{m}(s),z_{m}(s))|^{2}ds \nonumber\\
&\leq &
\ds\frac{4k_{2}^{2}}{\lambda_{1}^{1/2}}\ds\int_{0}^{T}|z_{m}(s)|^{2}\|w_{m}(s)\|^{2}
ds \nonumber\\
&\leq& \ds\frac{4k_{2}^{2}M_{1}}{\lambda_{1}^{1/2}}\ds\int_{0}^{T}|z_{m}(s)|^{2}ds
\nonumber\\
&= & \ds\frac{4k_{2}^{2}M_{1}}{\lambda_{1}^{1/2}}\|z_{m}\|^{2}_{L^{2}([0,T];H)}\leq
\ds\frac{4k_{2}^{2}M_{1}}{\lambda_{1}^{1/2}}M_{2}(T).\nonumber\\ \nonumber
\end{eqnarray}

To estimate the right-hand side of (\ref{estpraderivada}), we use the fact that $I_{h}(u)\in C([0,T];H)$ and so $A\mathcal{P}I_{h}(u)\in C([0,T];D(A)')$. Moreover, we have the two following estimates:
$$|I_{h}(w_{m})|\leq |I_{h}(w_{m})-w_{m}|+|w_{m}|\leq c_{1}h\|w_{m}\|+|w_{m}|,$$
$$\|AI_{h}(w_{m})\|_{D(A)'}=|A^{-1}AI_{h}(w_{m})|\leq |I_{h}w_{m}-w_{m}|+|w_{m}|\leq c_{1}h\|w_{m}\|+|w_{m}|.$$
Therefore, we conclude that
$$\left\|\ds\frac{dz_{m}}{dt}\ds\right\|^{2}_{L^{2}([0,T];H)}\leq M_{3}(\nu,\lambda_{1},\overline{f},
\alpha,T)\,\,\, \mbox{and}\,\,\,
\left\|\ds\frac{dw_{m}}{dt}\right\|^{2}_{L^{2}([0,T];D(A)')}
\leq M_{4}(\nu,\lambda_{1},\overline{f},\alpha,T),$$
for some $M_{3}$ and $M_{4}$.

Using the Aubin-Lions Compactness Theorem (see,e.g.,\cite{Constantin},\cite{Temam},\cite{Lion}) and the Banach-Alaoglu Theorem, we conclude that there exists a subsequence $\{w_{m_{j}}(t)\}_{j=1}^{\infty}$, that we denote with the same label $\{w_{m}(t)\}_{j=1}^{\infty}$
such that
\begin{equation}\label{conv1w}
w_{m}\rightarrow w\,\,\,\ \mbox{strongly in}\,\,\,\ L^{2}([0,T];V),
\end{equation}
\begin{equation}\label{conv2w}
w_{m}\rightarrow w\,\,\,\ \mbox{weakly in}\,\,\,\ L^{2}([0,T];D(A)),
\end{equation}
and equivalently,
\begin{equation}\label{conv1z}
z_{m}\rightarrow z\,\,\,\ \mbox{strongly in}\,\,\,\ L^{2}([0,T];V'),
\end{equation}
\begin{equation}\label{conv2z}
z_{m}\rightarrow z\,\,\,\ \mbox{weakly in}\,\,\,\ L^{2}([0,T];H).
\end{equation}
Using the same steps of Theorem \ref{teorema2002} from \cite{TitiHolmFoias}, one can show that the limit functions $w$ and $z$ satisfy (\ref{inversa}) and (\ref{espacoNS}).

We show next the continuous dependence on initial data of the solutions, and as a consequence, the uniqueness of the solution. Let $w$ and $\overline{w}$ be two solutions of $(\ref{versaocompridadaBtil})$ on the interval
$[0,T]$, with initial data $w(0)=w_{0}$ and $\overline{w}(0)=\overline{w}_{0}$, respectively. Denote $z=w+\alpha^{2}Aw$ and $\overline{z}=\overline{w}+\alpha^{2}A\overline{w}$. Denoting $\theta=\overline{w}-w$, we have
\begin{equation}\label{teta}
\frac{d}{dt}(\theta+\alpha^{2}A\theta)+\nu A(\theta+\alpha^{2}A\theta)+
\widetilde{B}(\overline{w},\overline{z})-\widetilde{B}(w,z)
=-\mu\mathcal{P}I_{h}(\theta)-\mu\alpha^{2}A\mathcal{P}(I_{h}(\theta)),
\end{equation}

Using that $\widetilde{B}(\overline{w},\overline{z})-\widetilde{B}(w,z)= \widetilde{B}(\theta,\overline{z})+\widetilde{B}(w,\theta+\alpha^{2}
A\theta)$, and the fact that (\ref{interdentrodo3}) holds in $L^{2}([0,T];D(A)')$, we take the inner product of (\ref{teta}) with $\theta$ and obtain
$$\begin{array}{rl}
\ds\frac{1}{2}\frac{d}{dt}(|\theta|^{2}
+\alpha^{2}\|\theta\|^{2})+\nu(\|\theta\|^{2}+\alpha^{2}|A\theta|^{2}) & \!\!\!\!\!
  +\langle\widetilde{B}(w,\theta+\alpha^{2}A\theta),
\theta\rangle_{D(A)',D(A)}\\
=& \!\!\! -\mu(I_{h}(\theta),\theta)-\mu\alpha^{2}(I_{h}(\theta),A\theta). \\
\end{array}$$

In the same way as was done in Theorem \ref{teoexismeu}, we estimate the following terms using Young's Inequality and the condition $2\mu c_{1}^{2}h^{2}<\nu$:
\begin{eqnarray}\label{unicaso1a}
-\mu(I_{h}(\theta),\theta) & = & \mu(\theta-I_{h}(\theta),\theta)-\mu|\theta|^{2}
\nonumber\\
&\leq & \mu|\theta-I_{h}(\theta)|\, |\theta|-\mu|\theta|^{2} \nonumber\\
&\leq &\mu c_{1}h\|\theta\|\, |\theta|-\mu|\theta|^{2}\\
&\leq & \ds\frac{\mu c_{1}^{2}h^{2}}{2}\|\theta\|^{2}+\ds\frac{\mu}{2}|\theta|^{2}-\mu|\theta|^{2}
 \nonumber\\
&\leq &\ds\frac{\nu}{4}\|\theta\|^{2}-\frac{\mu}{2}|\theta|^{2}.\nonumber
\nonumber
\end{eqnarray}
Similarly,
\begin{eqnarray}\label{unicaso1b}
-\mu\alpha^{2}(I_{h}\theta,A\theta)
& = & \mu\alpha^{2}(\theta-I_{h}\theta,A\theta)-\mu\alpha^{2}\|\theta\|^{2}
\nonumber\\
&\leq & \mu\alpha^{2}|\theta-I_{h}(\theta)|\, |A\theta|-\mu\alpha^{2}\|\theta\|^{2}
\nonumber\\
&\leq& \mu\alpha^{2}c_{1}^{2}h\|\theta\|\, |A\theta|-\mu\alpha^{2}\|\theta\|^{2} \\
&\leq& \alpha^{2}\mu c_{1}^{2}h^{2}|A\theta|^{2}+\ds\frac{\mu\alpha^{2}}{2}\|\theta\|^{2}-\mu\alpha^{2}\|\theta\|^{2}         \nonumber\\
&\leq&\ds\frac{\nu}{4}\alpha^{2}|A\theta|^{2}-\frac{\mu}{2}|\theta|^{2},\nonumber
\nonumber
\end{eqnarray}
where we used the assumption $2\mu c_{1}^{2}h^{2}\leq\ds\nu$.
Thus, from (\ref{teta})-(\ref{unicaso1b}),
\begin{equation}\label{eqparainspir}
\frac{1}{2}\frac{d}{dt}(|\theta|^{2}
+\alpha^{2}\|\theta\|^{2})+\frac{3\nu}{4}(\|\theta\|^{2}+\alpha^{2}|A\theta|)^{2}
  \leq|\langle\widetilde{B}(w,\theta+\alpha^{2}A\theta),
\theta\rangle_{V',V}|.
\end{equation}

To estimate the right-hand side of (\ref{eqparainspir}), we use (\ref{desiBtil8}), and Young's inequality to get

\begin{eqnarray}
|\langle\widetilde{B}(w,\theta+\alpha^{2}A\theta),
\theta\rangle_{V',V}| & \leq &   k_{3}\|w\|^{1/2}|Aw|^{1/2}|\theta|\,
\|\theta\|+k_{3}\alpha^{2}\|w\|^{1/2}|Aw|^{1/2}
|A\theta|\, \|\theta\| \nonumber\\
& + & k_{3}|Aw|\, |\theta|\, |\theta|^{1/2}\|\theta\|^{1/2}+k_{3}\alpha^{2}
|Aw|\, |A\theta|\, |\theta|^{1/2}\|\theta\|^{1/2} \nonumber\\
& \leq & \ds\frac{k_{3}^{2}}{\nu}\|w\|\, |Aw|\, |\theta|^{2}+\ds\frac{\nu}{4}\|\theta\|^{2}+\ds\frac{k_{3}^{2}}{\nu}\alpha^{2}\|w\|\, |Aw|\, \|\theta\|^{2}+\ds\frac{\nu}{4}\alpha^{2}|A\theta|^{2} \nonumber\\
 & + & \ds\frac{k_{3}^{2}}{\lambda_{1}^{1/2}\nu}|Aw|^{2}|\theta|^{2}+\ds\frac{\lambda_{1}^{1/2}\nu}{4}|\theta|\,
 \|\theta\|+\ds\frac{k_{3}^{2}}{\nu}\alpha^{2}|Aw|^{2}|\theta|\, \|\theta\|+\ds\frac{\nu}{4}\alpha^{2}|A\theta|^{2}. \nonumber
 \nonumber
\end{eqnarray}
Using Poincar{\'e}'s inequality, as well as Young's inequality again, it follows that
\begin{eqnarray}
|\langle\widetilde{B}(w,\theta+\alpha^{2}A\theta),
\theta\rangle_{V',V}| & \leq &  \ds\frac{k_{3}^{2}}{\lambda_{1}^{1/2}\nu}|Aw|^{2}|\theta|^{2}+\frac{k_{3}^{2}}{\lambda_{1}^{1/2}\nu}\alpha^{2}|Aw|^{2}\|\theta\|^{2}\nonumber\\
& + & \ds\frac{k_{3}^{2}}{\lambda_{1}^{1/2}\nu}|Aw|^{2}|\theta|^{2}+\ds\frac{k_{3}^{2}}{\lambda_{1}^{1/2}\nu}\alpha^{2}|Aw|^{2}\|\theta\|^{2}\nonumber\\
 & + &  \ds\frac{\nu}{2}(\|\theta\|^{2}+\alpha^{2}|A\theta|^{2}),\nonumber
 \nonumber
\end{eqnarray}
so we conclude the following estimate:
\begin{equation}\label{decima2}
|\langle\widetilde{B}(w,\theta+\alpha^{2}A\theta),
\theta\rangle_{V,V}|\leq \ds\frac{\nu}{2}(\|\theta\|^{2}+\alpha^{2}|A\theta|^{2})+
\ds\frac{2k_{3}^{2}}{\lambda_{1}^{1/2}\nu}|Aw|^{2}(|\theta|^{2}+\alpha^{2}\|\theta\|^{2}).
\end{equation}
Therefore, from above and (\ref{eqparainspir}) we have
\begin{eqnarray}
\ds\frac{1}{2}\frac{d}{dt}(
|\theta|^{2}+\alpha^{2}\|\theta\|^{2})+\!\!\!\!&&\!\!\!\!
 \!\!\!\! \ds\frac{\nu}{4}(\|\theta\|^{2}+ \alpha^{2}|A\theta|^{2})\nonumber\\
&\leq &
\ds\frac{2k_{3}^{2}}{\lambda_{1}^{1/2}\nu}|Aw|^{2}(|\theta|^{2}+\alpha^{2}\|\theta\|^{2}).\nonumber\\
 \nonumber
\end{eqnarray}
By Gronwall's inequality on the interval $[0,t]$, we obtain
$$|\theta(t)|^{2}+\alpha^{2}\|\theta(t)\|^{2}\leq(|\theta(0)|^{2}+\alpha^{2}
\|\theta(0)\|^{2})e^{\ds\frac{4k_{3}^{2}}{\lambda_{1}^{1/2}\nu}\ds\int_{0}^{t}|Aw(s)|^{2}ds}.$$

Since $w\in L^{2}([0,T];D(A))$, the above inequality implies the continuous dependence of the regular solution on the initial data, and in particular the uniqueness of a regular solution.

\Bo

For the case when $I_{h}$ is satisfying (\ref{become2}), we have the following result of well-posedness:
\begin{theorem}\label{teoexistmeu2}
Let $f\in H,w_{0}\in V$ and $\mu>0$ be given. Suppose that $I_{h}$ satisfies
(\ref{caso2}) (and hence $(\ref{become2})$) and $h$ is small enough such that the two conditions below are valid:
\begin{equation}
2\mu \overline{c}_{2}h^{2}<\nu
\end{equation}
and
\begin{equation}
2\mu c_{2}^{2}h^{4}<\nu\alpha^{2},
\end{equation}
where
$c_{2}>0$ is the constant given in $(\ref{become2})$ and $\overline{c}_{2}=\max\{c_{2},c_{2}^{2}\}$. Consider $u$ the solution of NS-$\alpha$ equations with initial data $u(0)=u_{0}\in V$, ensured by Theorem \ref{teorema2002}. Then the continuous data assimilation algorithm equations, $(\ref{eqcomInaexis})$, have a
 regular solution $w$ that satisfies
 \begin{equation}\label{espacoNS}
 w\in C([0,T];V)\cap L^{2}([0,T];D(A))\,\,\, \mbox{and}\,\,\, \ds\frac{dw}{dt}\in L^{2}([0,T];H),
 \end{equation}
 for any $T>0$.
\end{theorem}
{\it Proof.} The proof is similar to the proof of Theorem \ref{teoexismeu}. We start from inequality (\ref{iniciooutraexist}) by changing the way of using Young's inequality:
\begin{eqnarray}
\ds\frac{1}{2}\frac{d}{dt}(|w_{m}|^{2}+\alpha^{2}\|w_{m}\|^{2})
+ &&\!\!\!\!\!\!\!\!\!\!\!\! \nu(\|w_{m}\|^{2}+\alpha^{2}|Aw_{m}|^{2})
\leq\left(\ds\frac{1}{\mu}+\ds\frac{2\alpha^{2}}{\nu}\right)|\overline{f}|^{2} \nonumber \\
&+& \frac{\mu}{4}|w_{m}|^{2}+\ds\frac{\nu}{8}\alpha^{2}|Aw_{m}|^{2}\nonumber\\
& + & \ds\frac{\mu}{2}|w_{m}-I_{h}(w_{m})|^{2}+\ds\frac{\mu}{2}|w_{m}|^{2}-\mu|w_{m}|^{2}\nonumber\\
& + & \ds\frac{\mu^{2}\alpha^{2}}{\nu}|w_{m}-I_{h}(w_{m})|^{2}+\ds\frac{\nu}{4}\alpha^{2}|Aw_{m}|^{2}-\mu\alpha^{2}\|w_{m}\|^{2}.\nonumber
\nonumber
\end{eqnarray}

Applying $(\ref{become2})$, we obtain
\begin{eqnarray}
\ds\frac{1}{2}\frac{d}{dt}(|w_{m}|^{2}+\alpha^{2}\|w_{m}\|^{2})
+ &&\!\!\!\!\!\!\!\!\!\!\!\! \nu(\|w_{m}\|^{2}+\alpha^{2}|Aw_{m}|^{2})
\leq\left(\ds\frac{1}{\mu}+\ds\frac{2\alpha^{2}}{\nu}\right)|\overline{f}|^{2} \nonumber \\
&+& \frac{\mu}{4}|w_{m}|^{2}+\ds\frac{3\nu}{8}\alpha^{2}|Aw_{m}|^{2}\nonumber\\
& + & \ds\frac{\mu c_{2}^{2}h^{2}}{2}\|w_{m}\|^{2}+\ds\frac{\mu c_{2}^{2}h^{4}}{2}|Aw_{m}|^{2}-\ds\frac{\mu}{2}|w_{m}|^{2}\nonumber\\
& + & \ds\frac{\mu^{2}c_{2}^{2}h^{2}}{\nu}\alpha^{2}\|w_{m}\|^{2}+\ds\frac{\mu^{2}c_{2}^{2}h^{4}}{\nu}\alpha^{2}|Aw_{m}|^{2}-\mu\alpha^{2}\|w_{m}\|^{2}.\nonumber
\nonumber
\end{eqnarray}

The assumptions on $h$ that
$2\mu \overline{c}_{2}h^{2}<\nu$ and $2\mu c_{2}^{2}h^{4}<\nu\alpha^{2}$ imply that
$$2\mu c_{2}^{2}h^{2}\leq 2\mu \overline{c}_{2}\leq\nu, $$
and that
$$\mu^{2}c_{2}^{2}h^{4}=\mu c_{2}h^{2}\cdot\mu c_{2}h^{2}< \mu\overline{c}_{2}h^{2}\cdot\mu\overline{c}_{2}h^{2}<\frac{\nu}{2}\cdot\frac{\nu}{2}=\frac{\nu^{2}}{4}. $$

Therefore, the above implies
\begin{eqnarray}\label{paraintegrar3}
\ds\frac{1}{2}\frac{d}{dt}(|w_{m}|^{2}+\alpha^{2}\|w_{m}\|^{2})\!\!
& + &\!\! \ds\frac{\nu}{8}(\|w_{m}\|^{2}+\alpha^{2}|Aw_{m}|^{2})\nonumber\\
& \leq & \left(\ds\frac{1}{\mu}+\ds\frac{2\alpha^{2}}{\nu}\right)|\overline{f}|^{2}-\ds\frac{\mu}{4}(|w_{m}|^{2}+\alpha^{2}\|w_{m}\|^{2}).\nonumber
\end{eqnarray}
From here on, the proof follows similarly the steps of the proof of Theorem \ref{teoexismeu}, from inequality (\ref{paraintegrarnovo}) on.

\Bo

\section{The Convergence Theorems}
We derive now conditions on $\mu$ and $h$ in terms of physical parameters, such as $G$ and $\nu$, to guarantee the convergence, as $t\rightarrow\infty$, of the difference $w(t)-u(t)$ to zero, where $u$ is the solution of NS-$\alpha$ equations (ensured by Theorem $\ref{teorema2002}$) and $w$ solves (\ref{eqcomInaexis}).
 To do this, we first recall the following generalized Gronwall inequality, that can be found in \cite{JonesTiti}.

\begin{lemma}[Uniform Gronwall's Inequality]\label{gronwallgeneral}
Let $T>0$ be fixed, $\psi$ and $\beta$ be locally integrable real valued functions on $(0,\infty)$, satisfying the following conditions:
$$\liminf_{t\rightarrow\infty}\ds\int_{t}^{t+T}\beta(s)ds=\gamma>0 \,\,\, \mbox{and}
\,\,\, \limsup_{t\rightarrow\infty}\int_{t}^{t+T}\beta^{-}(s)ds=\Gamma<\infty,$$
where $\beta^{-}=\max\{-\beta,0\}$. Furthermore, assume that
$$\lim_{t\rightarrow\infty}\int_{t}^{t+T}\psi^{+}(s)ds=0,$$
where $\beta^{+}=\max\{\beta,0\}$. Suppose that $\xi$ is an absolutely continuous non-negative function on $(0,\infty)$, such that
$$\frac{d\xi(t)}{dt}+\beta(t)\xi(t)\leq\psi(t),$$
almost everywhere on $(0,\infty)$. Then $\xi(t)\rightarrow 0$ exponentially, as $t\rightarrow\infty$.
\end{lemma}

Our main result states

\begin{theorem}\label{teoremaconvergencia}
Let $u$ be a solution of the incompressible three-dimensional Navier-Stokes-$\alpha$ equations (\ref{versaocompridadaBtil}); and let $I_{h}:\dot{{H^{1}}}(\Omega)\rightarrow \dot{L^{2}}(\Omega)$ a linear map satisfying (\ref{become1}). Assume that $\mu$ is large enough satisfying
\begin{equation}\label{conddemu}
\mu>24k_{3}^{2}G^{2}+\frac{15c^{2}\nu G^{2}}{
\alpha^{2}},
\end{equation}
where $k_{3}$ is the constant appearing in $(\ref{desiBtil8})$. Moreover, assume that $h$ small enough such that
\begin{equation}\label{condproH}
h^{2}<\frac{\nu}{2\mu c_{1}^{2}}
\end{equation}
where $c_{1}$ is the constant given on (\ref{become1}). Then, the global unique solution $w$ of (\ref{eqcomInaexis}), given by Theorem \ref{teoexismeu}, satisfies $(w(t)-u(t))
\rightarrow 0$, as $t\rightarrow\infty$, in the $|\cdot|$ and $\|\cdot\|$-norms.
\end{theorem}

{\it Proof.} Considering $u$ the solution of the Navier-Stokes-$\alpha$ equations and denoting
 $\delta=w-u$, we have
 \begin{equation}\label{comdelta}
 \frac{d}{dt}(\delta+\alpha^{2}A\delta)+\nu A(\delta+\alpha^{2}A
\delta)+\widetilde{B}(w,z)-\widetilde{B}(u,v)=-\mu \mathcal{P}I_{h}(\delta)-\mu\alpha^{2}A\mathcal{P}I_{h}(\delta).
 \end{equation}
Note that
\begin{eqnarray}
\widetilde{B}(w,z)-\widetilde{B}(u,v) &= &
 \widetilde{B}(w,z-v)+\widetilde{B}(w-u,v)-\widetilde{B}(u,z-v)
+\widetilde{B}(u,z-v)\nonumber\\
& =& \widetilde{B}(\delta,\delta+\alpha^{2}A\delta)
+\widetilde{B}(\delta,v)+\widetilde{B}(u,\delta+\alpha^{2}A\delta).\nonumber
\nonumber
\end{eqnarray}

Taking the inner product of (\ref{comdelta}) with $\delta$ and using $(\ref{desiBtil4})$, we obtain
\begin{equation}\label{lionmagenes}
\begin{array}{rl}
\ds\frac{1}{2}\frac{d}{dt}(|\delta|^{2}+\alpha^{2}\|\delta\|^{2})
+\nu(\|\delta\|^{2}+
\alpha^{2}|A\delta|^{2})\,\ + \!\!\!&
  \langle\widetilde{B}(u,\delta+\alpha^{2}A\delta),
\delta\rangle_{D(A)',D(A)}  \\
&=-\mu(I_{h}(\delta),\delta)-\mu\alpha^{2}
\langle A\mathcal{P}I_{h}(\delta),\delta\rangle_{D(A)',D(A)}.\\
\end{array}
\end{equation}

Also, we have  $\langle A\mathcal{P}I_{h}(\delta),\delta\rangle_{D(A)',D(A)}=(I_{h}(\delta),A
\delta)$. Estimating the right-hand side terms of (\ref{lionmagenes}) using Young's inequality,
\begin{eqnarray}
-\mu(I_{h}(\delta),\delta) & = & \mu(\delta-I_{h}(\delta),\delta)-\mu|\delta|^{2}
\nonumber\\
&\leq & \mu|\delta-I_{h}(\delta)|\, |\delta|-\mu|\delta|^{2}\label{E2}\\
&\leq &\mu c_{1}h\|\delta\|\, |\delta|-\mu|\delta|^{2}\nonumber\\
&\leq & \frac{\mu}{2}|\delta|^{2}+\frac{\mu c_{1}^{2}h^{2}}{2}\|\delta\|^{2}-\mu|\delta|^{2}.\nonumber \nonumber
\end{eqnarray}

Similarly,
\begin{eqnarray}
-\mu\alpha^{2}(I_{h}\delta,A\delta)
& = & \mu\alpha^{2}(\delta-I_{h}(\delta),A\delta)-\mu\alpha^{2}\|\delta\|^{2}\nonumber\\
&\leq & \mu\alpha^{2}|\delta-I_{h}(\delta)|\, |A\delta|-\mu\alpha^{2}\|\delta\|^{2} \label{E3}\\
&\leq & \mu\alpha^{2}c_{1}h\|\delta\|\,|A\delta|-\mu\alpha^{2}\|\delta\|^{2}\nonumber\\
&=& \frac{\mu}{2}\alpha^{2}\|\delta\|^{2}+\frac{\mu c_{1}^{2}h^{2}}{2}\alpha^{2}\|\delta\|^{2}-\mu\alpha^{2}\|\delta\|^{2}.\nonumber
\nonumber
\end{eqnarray}

Therefore, from (\ref{lionmagenes}), (\ref{E2}) and (\ref{E3}) we have
\begin{eqnarray}
\ds\frac{1}{2}\frac{d}{dt}(|\delta|^{2}+\alpha^{2}\|\delta\|^{2})
& + & \nu(\|\delta\|^{2}+
\alpha^{2}|A\delta|^{2}) +
  \langle\widetilde{B}(u,\delta+\alpha^{2}A\delta),
\delta\rangle_{D(A)',D(A)}\nonumber\\
&\leq & \ds\frac{\mu}{2}|\delta|^{2}+\ds\frac{\mu c_{1}^{2}h^{2}}{2}\|\delta\|^{2}-\mu|\delta|^{2}
\label{duasestrelas}\\
& + & \ds\frac{\mu}{2}\alpha^{2}\|\delta\|^{2}+\ds\frac{\mu c_{1}^{2}h^{2}}{2}\alpha^{2}\|\delta\|^{2}-\mu\alpha^{2}\|\delta\|^{2}.\nonumber
\end{eqnarray}

The next step is to estimate the non-linear term. Using
$(\ref{desiBtil8})$ and Young's inequality, we have
\begin{eqnarray}
|\langle\widetilde{B}(u,\delta+\alpha^{2}A\delta),\delta\rangle_{V',V}| & \leq &
\frac{2k_{3}^{2}\lambda_{1}^{\frac{1}{2}}}{\nu}(\|u\|^{2}|\delta|^{2}+\alpha^{2}|Au|^{2}|\delta|^{2})
\nonumber\\
& + & \ds\frac{2k_{3}^{2}\lambda_{1}^{\frac{1}{2}}}{\nu}\alpha^{2}\|\delta\|^{2}\|u\|^{2}+\frac{\nu}{4}
(\|\delta\|^{2}+\alpha^{2}|A\delta|^{2})\label{E5}\\
& + & \ds\frac{5k_{3}^{2}}{2\nu\lambda_{1}^{1/2}}|Au|^{2}
(|\delta|^{2}+\alpha^{2}\|\delta\|^{2}).\nonumber\\
\nonumber
\end{eqnarray}
Thus,
\begin{eqnarray}
\ds\frac{1}{2}\frac{d}{dt}(|\delta|^{2}+\alpha^{2}\|\delta\|^{2})
&\!\!\! + &\!\!\! \nu(\|\delta\|^{2}+
\alpha^{2}|A\delta|^{2})\nonumber \\
& \leq & \!\!\! \frac{4k_{3}^{2}\lambda_{1}^{\frac{1}{2}}}{\nu}
(\|u\|^{2}+\alpha^{2}|Au|^{2})(|\delta|^{2}+\alpha^{2}\|\delta\|^{2})\nonumber\\
& + & \!\!\! \frac{5k_{3}^{2}}{2\nu\lambda_{1}^{1/2}}|Au|^{2}
(|\delta|^{2}+\alpha^{2}\|\delta\|^{2})+\frac{\nu}{4}
(\|\delta\|^{2}+\alpha^{2}|A\delta|^{2})\nonumber\\
& + & \!\!\! \ds\frac{\mu}{2}|\delta|^{2}+\ds\frac{\mu c_{1}^{2}h^{2}}{2}\|\delta\|^{2}-\mu|\delta|^{2}
+ \ds\frac{\mu}{2}\alpha^{2}\|\delta\|^{2}+\ds\frac{\mu c_{1}^{2}h^{2}}{2}\alpha^{2}\|\delta\|^{2}-\mu\alpha^{2}\|\delta\|^{2}.\nonumber
\nonumber
\end{eqnarray}

Since by assumption $2\mu c_{1}^{2}h^{2}<\ds\nu$, we have
\begin{eqnarray}
\ds\frac{1}{2}\frac{d}{dt}(|\delta|^{2}+\alpha^{2}\|\delta\|^{2})
&\!\!\! + &\!\!\! \frac{\nu}{2}(\|\delta\|^{2}+
\alpha^{2}|A\delta|^{2})\nonumber \\
& \leq & \!\!\! \frac{4k_{3}^{2}\lambda_{1}^{\frac{1}{2}}}{\nu}
(\|u\|^{2}+\alpha^{2}|Au|^{2})(|\delta|^{2}+\alpha^{2}\|\delta\|^{2})\label{usarprox}\\
& + & \!\!\! \frac{5k_{3}^{2}}{2\nu\lambda_{1}^{1/2}}|Au|^{2}
(|\delta|^{2}+\alpha^{2}\|\delta\|^{2})
-\frac{\mu}{2}(|\delta|^{2}+\alpha^{2}\|\delta\|^{2}).\nonumber\\
\end{eqnarray}
Therefore, we conclude that
\begin{equation}\label{edogronwall}
\ds\frac{d}{dt}(|\delta(t)|^{2}+\alpha^{2}\|\delta(t)\|^{2})+\beta(t)(|\delta(t)|^{2}+\alpha^{2}\|\delta(t)\|^{2})\leq 0,
\end{equation}
where
$$\beta(t)=\mu-\frac{8k_{3}^{2}\lambda_{1}^{\frac{1}{2}}}{\nu}(\|u(t)\|^{2}+\alpha^{2}|Au(t)|^{2})-
\frac{5k_{3}^{2}}{\nu\lambda_{1}^{1/2}}|Au(t)|^{2}.$$

To make use of Lemma \ref{gronwallgeneral}, we note that for $T>0$,
\begin{equation}\label{betadesig}
\ds\int_{t}^{t+T}\beta(s)ds=\mu T-\frac{8c^{2}\lambda_{1}^{\frac{1}{2}}}{\nu}\ds\int_{t}^{t+T}
\|u(s)\|^{2}+\alpha^{2}|Au(s)|^{2}ds-\frac{5c^{2}}{\nu
\lambda_{1}^{1/2}}\cdot\frac{1}{\alpha^{2}}\ds\int_{t}^{t+T}\alpha^{2}|Au(s)|^{2}ds.
\end{equation}

Taking $T=\ds\frac{1}{\nu\lambda_{1}}$ in Proposition \ref{estdosatratores}, we have for $t\geq t_{0}$ (where $t_{0}$ is given in Proposition \ref{estdosatratores}),
\begin{equation}\label{duasflores}
-\frac{5c^{2}}{\nu
\lambda_{1}^{1/2}}\cdot\frac{1}{\alpha^{2}}\ds\int_{t}^{t+T}\alpha^{2}|Au(s)|^{2}ds\geq -
\frac{15k_{3}^{2}G^{2}}{\lambda_{1}\alpha^{2}}.
\end{equation}

Therefore, to guarantee that $\ds\liminf_{t\rightarrow\infty}\ds\int_{t}^{t+T}\beta(s)ds=\gamma>0$, it is sufficient to require
$$\ds\frac{\mu}{\nu\lambda_{1}}-24k_{3}^{2}G^{2}-
\frac{15k_{3}^{2}G^{2}}{\lambda_{1}\alpha^{2}}>0,$$
 which is given by assumption (\ref{conddemu}).

 Finally, taking $\psi\equiv 0$ in Lemma \ref{gronwallgeneral}, we conclude that
$$|\delta(t)|^{2}+\alpha^{2}\|\delta(t)\|^{2}\longrightarrow 0, \,\,\, \mbox{as} \,\,\, t\rightarrow\infty,$$
i.e., $(w(t)-u(t))\rightarrow 0$, in $L^{2}$ and $H^{1}$-norms, exponentially in time, as $t\rightarrow\infty$.

\Bo

Note that Theorem $\ref{teoremaconvergencia}$ is valid for any arbitrary initial value $w_{0}$ of $w$, and this is the main advantage of this type of assimilation, because it overcomes the difficulty coming from the lack of information on the initial data of the reference solution.

We consider now the case where the interpolant operator $I_{h}$ satisfies the approximation property (\ref{become2}). Then we have the following theorem of convergence:
\begin{theorem}\label{teoconvergencia2}
Let $u$ be the solution of Navier-Stokes-$\alpha$ equations (\ref{versaocompridadaBtil}) and $I_{h}:\dot{H^{2}}(\Omega)\rightarrow \dot{L^{2}} (\Omega)$ be a linear interpolant operator satisfying (\ref{become2}). Suppose that $\mu$ is large enough satisfying
\begin{equation}\label{condteo2}
\mu>24k_{3}^{2}\nu\lambda_{1}G^{2}+
\ds\frac{15k_{3}^{2}\nu G^{2}}{\alpha^{2}},
\end{equation}
 and $h$ small enough such that
\begin{equation}\label{condidupla}
2\mu \overline{c}_{2}h^{2}\leq\nu\,\,\, \mbox{and}\,\,\, 2\mu c_{2}^{2}h^{4}\leq\nu\alpha^{2},
\end{equation}
where $\overline{c}_{2}=\max\{c_{2},c_{2}^{2}\}$.
 Then, the global unique solution $w$ of (\ref{eqcomInaexis}), given by Theorem \ref{teoexistmeu2}, satisfies $(w(t)-u(t))
\rightarrow 0$, as $t\rightarrow\infty$, in the $L^{2}$ and $H^{1}$-norms.
\end{theorem}
{\it Proof.} The idea of the proof is the same as in Theorem $\ref{teoremaconvergencia}$, except for the fact that we need here to estimate $-\mu(I_{h}\delta,\delta)$ and
$-\mu\alpha^{2}(I_{h}\delta,A\delta)$ as follows:
\begin{eqnarray}
-\mu(I_{h}\delta,\delta) & = & \mu(\delta-I_{h}\delta,\delta)-\mu|\delta|^{2}
\nonumber\\
&\leq& \mu|\delta-I_{h}(\delta)||\delta|-\mu|\delta|^{2}\nonumber\\
&\leq & \ds\frac{\mu}{2}|\delta-I_{h}(\delta)|^{2}+\ds\frac{\mu}{2}|\delta|^{2}-\mu|\delta|^{2}\\
&\leq & \frac{\mu c_{2}^{2}h^{2}}{2}\|\delta\|^{2}+\frac{\mu}{2}c_{2}^{2}h^{4}|A\delta|^{2}-\frac{\mu}{2}|\delta|^{2}.
\nonumber\\ \nonumber
\end{eqnarray}

By the assumption (\ref{condidupla}) we have $\mu c_{2}^{2}h^{2}\leq\mu \overline{c}_{2}h^{2}<\ds\frac{\nu}{2}$ and $\mu c_{2}^{2}h^{4}<\ds\frac{\nu\alpha^{2}}{2}$. Therefore, (\ref{condidupla}) implies
\begin{equation}\label{menosmu1}
-\mu(I_{h}\delta,\delta)\leq\ds\frac{\nu}{4}\|\delta\|^{2}+\ds\frac{\nu}{4}\alpha^{2}|A\delta|^{2}-\ds\frac{\mu}{2}|\delta|^{2}=\frac{\nu}{4}(\|\delta\|^{2}+\alpha^{2}|A\delta|^{2})-\ds\frac{\mu}{2}|\delta|^{2}.
\end{equation}

Furthermore, by Cauchy-Schwarz and Young inequalities, we have
\begin{eqnarray}
-\mu\alpha^{2}(I_{h}\delta,A\delta)
& = & \mu\alpha^{2}(\delta-I_{h}\delta,A\delta)-\mu\alpha^{2}\|\delta\|^{2}\nonumber\\
&\leq & \mu\alpha^{2}|\delta-I_{h}(\delta)||A\delta|-\mu\alpha^{2}\|\delta\|^{2}\nonumber\\
&\leq & \ds\frac{\mu^{2}\alpha^{2}}{\nu}|\delta-I_{h}(\delta)|^{2}+\ds\frac{\nu}{4}\alpha^{2}|A\delta|^{2}-\mu\alpha^{2}\|\delta\|^{2}
\nonumber\\
&\leq &\ds\frac{\mu^{2}\alpha^{2}c_{2}^{2}h^{2}}{\nu}\|\delta\|^{2}+\frac{\mu^{2}\alpha^{2}c_{2}^{2}h^{4}}{\nu}|A\delta|^{2}+\ds\frac{\nu}{4}\alpha^{2}|A\delta|^{2}-\mu\alpha^{2}\|\delta\|^{2}.\nonumber
\nonumber
\end{eqnarray}

Thanks to (\ref{condidupla}) we have, as a result, that $\mu c_{2}^{2} h^{2}<\ds\frac{\nu}{2}$, and hence
$$\mu^{2}c_{2}^{2}h^{2}=\mu\cdot \mu c_{2}^{2}h^{2}\leq \frac{\mu\nu}{2}.$$

Consequently,
$$\mu^{2}c_{2}^{2}h^{4}=\mu c_{2}h^{2}\cdot \mu c_{2}h^{2}\leq \mu\overline{c}_{2}h^{2}\cdot \mu\overline{c}_{2}h^{2}<\ds\frac{\nu}{2}\cdot\ds\frac{\nu}{2}=\ds\frac{\nu^{2}}{4}.$$

Therefore,
\begin{eqnarray}\label{menosmu2}
-\mu\alpha^{2}(I_{h}\delta,A\delta) & \leq &
\ds\frac{\mu\nu\alpha^{2}}{2\nu}\|\delta\|^{2}+\ds\frac{\nu^{2}\alpha^{2}}{4\nu}|A\delta|^{2}+\ds\frac{\nu}{4}\alpha^{2}|A\delta|^{2}-\mu\alpha^{2}\|\delta\|^{2}\nonumber\\
& = & \ds\frac{\mu}{2}\alpha^{2}\|\delta\|^{2}+\frac{\nu}{2}\alpha^{2}|A\delta|^{2}-\mu\alpha^{2}\|\delta\|^{2}\\
& = & \ds\frac{\nu}{2}\alpha^{2}|A\delta|^{2}-\ds\frac{\mu}{2}\alpha^{2}\|\delta\|^{2}.\nonumber
\nonumber
\end{eqnarray}

Using (\ref{E5}), together with $(\ref{menosmu1})$ and $(\ref{menosmu2})$ we obtain
\begin{eqnarray}
\ds\frac{1}{2}\frac{d}{dt}(|\delta|^{2}+\alpha^{2}\|\delta\|^{2})
&+ &\!\!\! \nu(\|\delta\|^{2}+
\alpha^{2}|A\delta|^{2})\nonumber \\
& \leq & \!\!\! \frac{4c^{2}\lambda_{1}^{\frac{1}{2}}}{\nu}
(\|u\|^{2}+\alpha^{2}|Au|^{2})(|\delta|^{2}+\alpha^{2}\|\delta\|^{2})\nonumber\\
& + & \!\!\! \frac{5c^{2}}{2\nu\lambda_{1}^{1/2}}|Au|^{2}
(|\delta|^{2}+\alpha^{2}\|\delta\|^{2})+\frac{\nu}{4}
(\|\delta\|^{2}+\alpha^{2}|A\delta|^{2}) \\
& + & \frac{\nu}{4}(\|\delta\|^{2}+\alpha^{2}|A\delta|^{2})-\ds\frac{\mu}{2}|\delta|^{2}+\ds\frac{\nu}{2}\alpha^{2}|A\delta|^{2}-\ds\frac{\mu}{2}\alpha^{2}\|\delta\|^{2}.\nonumber
\nonumber
\end{eqnarray}
As a result, it follows that
\begin{eqnarray}
\ds\frac{d}{dt}(|\delta|^{2}+\alpha^{2}\|\delta\|^{2})
& \leq & \!\!\! \frac{8c^{2}\lambda_{1}^{\frac{1}{2}}}{\nu}
(\|u\|^{2}+\alpha^{2}|Au|^{2})(|\delta|^{2}+\alpha^{2}\|\delta\|^{2})\nonumber\\
& + & \!\!\! \frac{5c^{2}}{2\nu\lambda_{1}^{1/2}}|Au|^{2}
(|\delta|^{2}+\alpha^{2}\|\delta\|^{2})-
\mu(|\delta|^{2}+\alpha^{2}\|\delta\|^{2}).
\nonumber
\end{eqnarray}

As in Theorem \ref{teoremaconvergencia}, we have
\begin{equation}
\ds\frac{d}{dt}(|\delta(t)|^{2}+\alpha^{2}\|\delta(t)\|^{2})+\beta(t)(|\delta(t)|^{2}+\alpha^{2}\|\delta(t)\|^{2})\leq 0,
\end{equation}
where
$$\beta(t)=\mu-\frac{8c^{2}\lambda_{1}^{\frac{1}{2}}}{\nu}(\|u(t)\|^{2}+\alpha^{2}|Au(t)|^{2})-
\frac{5c^{2}}{\nu\lambda_{1}^{1/2}}|Au(t)|^{2}.$$

 Using exactly the same calculations as in Theorem \ref{teoremaconvergencia}, we make use of Lemma \ref{gronwallgeneral} to conclude that for $\mu$ large enough satisfying (\ref{condteo2}) and $h$ is small enough such that (\ref{condidupla}) holds, then
$$|\delta(t)|^{2}+\alpha^{2}\|\delta(t)\|^{2}\longrightarrow 0, \,\,\, \mbox{as} \,\,\, t\rightarrow\infty,$$
 exponentially in time, which is the desired conclusion.

 \Bo

\section{Examples of interpolant operators}

In this section we give some examples of interpolant operators satisfying the approximating identity inequalities (\ref{caso1}) or (\ref{caso2}). For two-dimensional cases, similar examples that will be considered here are found in \cite{AzouaniOlsonTiti}

It is a simple exercise of Fourier Analysis to prove that the interpolant $I_{h}:\dot{H^{1}}(\Omega)\rightarrow \dot{L^{2}}(\Omega)$ given by the projection onto the low Fourier modes with wave numbers $k$ such that $|k|\leq\lfloor\lambda_{1}^{-\frac{1}{2}}h^{-1}\rfloor$:
$$I_{h}\varphi=P_{k}\varphi=\ds\sum_{|k|\leq \lfloor\lambda_{1}^{-\frac{1}{2}}h^{-1}\rfloor}\widehat{\varphi}_{k}\phi_{k},  \,\,\,\, \mbox{where}\,\,\,\,
\varphi(x)=\ds\sum_{k\in\mathbb{Z}^{3}\backslash\{0\}}
\widehat{\varphi}_{k}\phi_{k}(x),$$
satisfies
\begin{equation}\label{examples}
\|\varphi-I_{h}\varphi\|^{2}_{L^{2}(\Omega)}\leq c_{1}^{2}h^{2}\|\nabla\varphi\|^{2}_{L^{2}(\Omega)}.
\end{equation}

Another example of an interpolant that satisfies the condition (\ref{examples}) is that given by volume elements, which is physically important and was studied in \cite{JonesTiti} in the context of 2D Navier-Stokes equations. We divide the periodic domain $\Omega=[0,L]^{3}$ in $\Omega_{k},k=1,...,N$, where $\Omega_{k}$ is the cube with edge $\ds\frac{L}{\sqrt[3]{N}}$, and so $|\Omega_{k}|=\ds\frac{L^{3}}{N}\ds$.
 Recalling that the local average of $u$ in $\Omega_{k}$ is defined as
$\langle u\rangle_{\Omega_{k}}=\frac{1}{|\Omega_{k}|}\int_{\Omega_{k}}u(x)dx$, we construct $I_{h}$ as follows:
\begin{equation}\label{media}
I_{h}(\varphi(x))=\ds\sum_{k=1}^{N}\langle \varphi\rangle_{\Omega_{k}}\chi_{\Omega_{k}}(x),
\end{equation}
 where $h=\ds\frac{L}{\sqrt[3]{N}}$. Here we suppose that the average values of $\varphi$ on each of the $\Omega_{k}$'s is given. To prove that $I_{h}$ satisfies (\ref{examples}), we generalize the result obtained in \cite{JonesTiti} for the two-dimensional case $\Omega=[0,L]^{2}$: for all $u\in H^{1}([0,L]^{2})$,
  $$\|u\|^{2}_{L^{2}(K_{j})}\leq l^{2}\langle u\rangle^{2}_{K_{j}}+\frac{l^{2}}{3}
\|\nabla u\|_{L^{2}(K_{j})}^{2}.$$
 where the domain $[0,L]^{2}$ has been divided into $N$ squares $K_{j}$ with side $l=\frac{L}{\sqrt[2]{N}},j=1,...,N$. Generalizing, we can prove that, in 3D,
 \begin{equation}\label{3dparavolume}
 \|u\|^{2}_{L^{2}(\Omega_{j})}\leq l^{3}\langle u\rangle_{\Omega_{j}}^{2}+\frac{l^{2}}{3}\|\nabla u\|^{2}
_{L^{2}(\Omega_{j})},
 \end{equation}
for all $j=1,...,N$. Using (\ref{3dparavolume}), we prove that
$$\|\varphi-I_{h}\varphi\|^{2}_{L^{2}(\Omega)}\leq\ds\frac{1}{3}h^{2}\|\nabla\varphi\|^{2}_{L^{2}(\Omega)}. $$

To conclude, we examine the most physically interesting example of an interpolant $I_{h}$ which satisfies
\begin{equation}\label{novadesint}
 \|\varphi-I_{h}\varphi\|^{2}_{L^{2}(\Omega)}\leq c_{2}^{2}h^{2}\|\nabla\varphi\|^{2}_{L^{2}(\Omega)}+c_{2}^{2}h^{4}\|\varphi\|^{2}_{H^{2}(\Omega)},
 \end{equation}
namely the interpolant obtained using measurements at a discrete set of nodal points in $\Omega=[0,L]^{3}$.

Indeed, similarly to what was done in the previous example of volume elements, to construct such an interpolant using nodal values we divide the domain $\Omega$ in $N$ cubes of edge $\ds\frac{L}{\sqrt[3]{N}}$, for $N\in\mathbb{N}$ and thus $|\Omega_{j}|=\ds\frac{L^{3}}{N},j=1
,...,N$, where $\Omega_{j}$ denote the $j-$th cube and $\Omega=\ds\cup_{j=1}^{N}
\Omega_{j}$. We then regard arbitrary points $x_{j}\in \Omega_{j}$ that represent the points where observational measurements of the velocity of the flow are done.

Define this interpolant as
\begin{equation}\label{nodal}
I_{h}\varphi(x)=\ds\sum_{k=1}^{N}\varphi(x_{k})\chi_{\Omega_{k}}(x).
\end{equation}

To prove that the interpolant above satisfies (\ref{novadesint}), we make use of the following two lemmas:

\begin{lemma}\label{lema6.1}
Let $\overline{Q}=[0,\Lambda]\times[0,d]$ and $u\in
H^{1}(\overline{Q})$. Then
\begin{equation}
\ds\int_{0}^{\Lambda}|u(x,0)|^{2}dx\leq\ds\frac{2}{d}\|u\|^{2}_{L^{2}(\overline{Q})}+d\left\|\ds\frac{\partial u}{\partial y}\right\|^{2}_{L^{2}(\overline{Q})},
\end{equation}
and similarly,
\begin{equation}
\ds\int_{0}^{d}|u(0,y)|^{2}dy\leq\ds\frac{2}{\Lambda}\|u\|^{2}_{L^{2}(\overline{Q})}+\Lambda\left\|\ds\frac{\partial u}{\partial x}\right\|^{2}_{L^{2}(\overline{Q})}.
\end{equation}
\end{lemma}
{\bf Proof:} This can be found in Lemma 6.1 of \cite{JonesTiti2}.

\begin{lemma}\label{difdospontos}
Let $\overline{\Omega}=[0,l]\times[0,l]\times[0,l]$ and $x$ and $z$ be two points of $\overline{\Omega}$, where the third coordinates of $x$ and $z$ are the same, ie, $x=(x_{1},y_{1},z_{1})$ and $z=(x_{2},y_{2},z_{1})$. Then for every $\varphi\in H^{2}(\overline{\Omega})$, we have
\begin{equation}\label{terccoordigual}
|\varphi(x)-\varphi(z)|\leq\ds\frac{2}{l^{1/2}}\left(
4\|\nabla\varphi\|^{2}_{L^{2}(\overline{\Omega})}+l^{2}
\left\|\frac{\partial^{2}\varphi}{\partial x\partial y}\right\|^{2}_{L^{2}(\overline{\Omega})}\right)^{\frac{1}{2}}.
\end{equation}

Similarly, if $y$ and $z$ are two points in $\overline{\Omega}$ such that the second coordinate of $y$ and $z$ are the same, ie, $y=(x_{2},y_{2},z_{2})$ and $z=(x_{3},y_{2},z_{3})$, then
\begin{equation}\label{segundcoordigual}
|\varphi(y)-\varphi(z)|\leq\ds\frac{2}{l^{1/2}}\left(
4\|\nabla\varphi\|^{2}_{L^{2}(\overline{\Omega})}+l^{2}
\left\|\frac{\partial^{2}\varphi}{\partial x\partial z}\right\|^{2}_{L^{2}(\overline{\Omega})}\right)^{\frac{1}{2}},
\end{equation}

for every $\varphi\in H^{2}(\overline{\Omega})$.
\end{lemma}

{\it Proof.}
We will show only the first estimate, and the second one is analogous. We begin by considering the square $Q=[0,l]\times[0,l]$. For any two points in $\overline{\Omega}$
 of the form
$(x_{1},y,z_{1})$ and $(x_{2},y,z_{1})$, with $y\in[0,l]$, we have
\begin{equation}\label{difcomzfixo}
|\varphi(x_{1},y,z_{1})  - \varphi(x_{2},y,z_{1})|^{2} =
\left|\ds\int_{x_{1}}^{x_{2}}\frac{\partial\varphi}{\partial x}(s,y,z_{1})ds\right|^{2}
\leq l\left\|\frac{\partial\varphi}{\partial x}(\cdot,y,z_{1})\right\|^{2}_{L^{2}([0,l])}
\end{equation}

Since the third coordinate $z_{1}$ is fixed and the points $(x_{1},y,z_{1})$ and $\varphi(x_{2},y,z_{1})$ are in a plane parallel to the $xy$ plane, we can apply Lemma
 \ref{lema6.1} for $\ds\frac{\partial\varphi}{\partial x}(\cdot,y,z_{1})$, with $d$ replaced with the maximal distance of the $y$-coordinate of the points $(x_{1},y,z_{1}),(x_{2},y,z_{1})$ from the horizontal walls; ie,
 $$l\geq d=\max\{y,l-y\}\geq\frac{l}{2}$$
and therefore
$$\ds\int_{0}^{l}\left|\frac{\partial\varphi}{\partial x}(x,y,z_{1})\right|^{2}dx\leq\ds\frac{4}{l}\left\|\frac{\partial\varphi}{\partial x}\right\|^{2}_{L^{2}(Q)}+l\left\|\ds\frac{\partial^{2}u}{\partial y\partial x}\right\|^{2}_{L^{2}(Q)},$$
since $\ds\frac{1}{d}\leq\frac{2}{l}$. Then we have
\begin{equation}\label{comzfixo}
l\left\|\frac{\partial\varphi}{\partial x}(\cdot,y,z_{1})\right\|^{2}_{L^{2}([0,l])}\leq 4\left\|\frac{\partial\varphi}{\partial x}\right\|^{2}_{L^{2}(Q)}+l^{2}\left\|\ds\frac{\partial^{2}\varphi}{\partial y\partial x}\right\|^{2}_{L^{2}(Q)}.
\end{equation}

Plugging $(\ref{comzfixo})$ into $(\ref{difcomzfixo})$, we have that
\begin{equation}
|\varphi(x_{1},y,z_{1})-\varphi(x_{2},y,z_{1})|^{2}\leq
4\left\|\ds\frac{\partial\varphi}{\partial x} \right\|^{2}_{L^{2}(Q)}+l^{2}\left\|\ds\frac{\partial^{2}\varphi}{\partial y\partial x} \right\|^{2}_{L^{2}(Q)}.
\end{equation}

By symmetry, we have the similar inequality for points of the form $(x,y_{1},z_{1})$ and $(x,y_{2},z_{1})$, where $x\in(0,l)$:
\begin{equation}
|\varphi(x,y_{1},z_{1})-\varphi(x,y_{2},z_{1})|^{2}\leq
4\left\|\ds\frac{\partial\varphi}{\partial y} \right\|^{2}_{L^{2}(Q)}+l^{2}\left\|\ds\frac{\partial^{2}\varphi}{\partial x\partial y} \right\|^{2}_{L^{2}(Q)}.
\end{equation}

Thus
$$\begin{array}{lll}
|\varphi(x)-\varphi(z)|^{2} & = & |\varphi(x_{1},y_{1},z_{1})-\varphi(x_{2},y_{2},z_{1})|^{2}\nonumber\\
\\
&\leq &(|\varphi(x_{1},y_{1},z_{1})-\varphi(x_{2},y_{1},z_{1})|+|\varphi(x_{2},y_{1},z_{1})-\varphi(x_{2},y_{2},z_{1})|)^{2}\nonumber\\
\\
& \leq & 2|\varphi(x_{1},y_{1},z_{1})-\varphi(x_{2},y_{1},z_{1})|^{2}+2|\varphi(x_{2},y_{1},z_{1})-\varphi(x_{2},y_{2},z_{1})|^{2}\nonumber\\
\\
& \leq & 2\left(4\left\|\ds\frac{\partial\varphi}{\partial x} \right\|^{2}_{L^{2}(Q)}+l^{2}\left\|\ds\frac{\partial^{2}\varphi}{\partial y\partial x} \right\|^{2}_{L^{2}(Q)}\right)+2\left(4\left\|\ds\frac{\partial\varphi}{\partial y} \right\|^{2}_{L^{2}(Q)}+l^{2}\left\|\ds\frac{\partial^{2}\varphi}{\partial x\partial y} \right\|^{2}_{L^{2}(Q)}\right),\nonumber\\
\end{array}$$
and it follows that
\begin{equation}\label{difpontos2d}
|\varphi(x)-\varphi(z)|^{2}\leq 4\left(4\|\nabla\varphi\|^{2}_{L^{2}(Q)}+l^{2}\left\|\ds\frac{\partial^{2}\varphi}{\partial y\partial x} \right\|^{2}_{L^{2}(Q)}\right).
\end{equation}

Our aim is to obtain estimates in $L^{2}(\overline{\Omega})$, where $\overline{\Omega}=[0,l]\times[0,l]\times[0,l]$ instead of $L^{2}(Q)$. For this, we integrate $(\ref{difpontos2d})$ from 0 to $l$ in $z$-coordinate:
\begin{eqnarray}
\int_{0}^{l}|\varphi(x_{1},y_{1},z_{1})-\varphi(x_{2},y_{2},z_{1})|dz & \leq & 2\int_{0}^{l}\left(
4\|\nabla\varphi(\cdot,\cdot,z)\|^{2}_{L^{2}(Q)}+l^{2}\left\|\ds\frac{\partial^{2}\varphi}{\partial y\partial x}(\cdot,\cdot,z) \right\|^{2}_{L^{2}(Q)}\right)^{\frac{1}{2}}dz\nonumber\\
& \leq & 2\int_{0}^{l}\left(
4\|\nabla\varphi(\cdot,\cdot,z)\|^{2}_{L^{2}(Q)}+l^{2}\left\|\ds\frac{\partial^{2}\varphi}{\partial y\partial x}(\cdot,\cdot,z) \right\|^{2}_{L^{2}(Q)}dz\right)^{\frac{1}{2}} l^{\frac{1}{2}}.\nonumber\\
\nonumber
\end{eqnarray}

Therefore,
$$l|\varphi(x_{1},y_{1},z_{1})-\varphi(x_{2},y_{2},z_{1})|\leq 2l^{\frac{1}{2}}\left(4
\|\nabla\varphi\|^{2}_{L^{2}(\overline{\Omega})}+l^{2}\left\|\frac{\partial^{2}\varphi}{\partial x\partial y}\right\|^{2}_{L^{2}(\overline{\Omega})}\right)^{\frac{1}{2}},$$
ie,
\begin{equation}
|\varphi(x_{1},y_{1},z_{1})-\varphi(x_{2},y_{2},z_{1})|^{2}\leq \frac{4}{l}\left(4
\|\nabla\varphi\|^{2}_{L^{2}(\overline{\Omega})}+l^{2}\left\|\frac{\partial^{2}\varphi}{\partial x\partial y}\right\|^{2}_{L^{2}(\overline{\Omega})}\right).
\end{equation}
and we have the desired conclusion.

\Bo

We are ready to prove that the interpolant $I_{h}$ constructed using measurements at nodal points satisfies
$(\ref{novadesint})$:

\begin{lemma}
For all $\varphi\in D(A)$, the interpolant $I_{h}$ defined in $(\ref{nodal})$ satisfies
\begin{equation}\label{intdentrolema}
|\varphi-I_{h}\varphi|^{2}\leq 32h^{2}\|\varphi\|^{2}+4h^{4}|A\varphi|^{2}
\end{equation}
where $h=L/ \sqrt[3]{N}$.
\end{lemma}
{\it Proof.} Note that
\begin{eqnarray}\label{contagigante}
|\varphi & \!\!\! - & \!\!\!\ds\sum_{k=1}^{N}\varphi(x_{k})\chi_{\Omega_{k}}|^{2} =
\ds\int_{\Omega}|\varphi(x)-\varphi(x_{k})\chi_{\Omega_{k}}(x)|^{2}dx
\nonumber\\
& =\!\! & \ds\int_{\Omega}|\varphi(x)\displaystyle\sum_{k=1}^{N}
 \chi_{\Omega_{k}}(x)-\sum_{k=1}^{N}
 \varphi(x_{k})\chi_{\Omega_{k}}(x)|^{2}dx\nonumber
\nonumber\\
 & =\!\! & \ds\int_{\Omega}\left[\sum_{k=1}^{N}\varphi(x)
 \chi_{\Omega_{k}}(x)-\sum_{k=1}^{N}
\varphi(x_{k})\chi_{\Omega_{k}}(x)\right]
\left[\sum_{j=1}^{N}\varphi(x)
 \chi_{\Omega_{j}}(x)-\sum_{j=1}^{N}
\varphi(x_{k})\chi_{\Omega_{k}}(x)\right]dx.\nonumber\\
\nonumber
\end{eqnarray}

Since $\chi_{\Omega_{k}}(x)\chi_{\Omega_{j}}(x)=\chi_{\Omega_{k}}(x)\delta_{kj}$, we have
\begin{eqnarray}\label{somatorianodal1}
|\varphi & \!\!\!\! - & \!\!\!\!\ds\sum_{k=1}^{N}\varphi(x_{k})\chi_{\Omega_{k}}|^{2}\leq
\ds\int_{\Omega}\left[\sum_{k=1}^{N}(\varphi(x)-\varphi(x_{k}))
\chi_{\Omega_{k}}(x)\right]\left[\sum_{j=1}^{N}(\varphi(x)-\varphi(x_{k}))
\chi_{\Omega_{j}}(x)\right]dx\nonumber\\
& =\!\! & \ds\int_{\Omega}\sum_{k=1}^{N}\sum_{j=1}^{N}(\varphi(x)-
\varphi(x_{k}))(\varphi(x)-
\varphi(x_{k}))\chi_{\Omega_{k}}(x)\chi_{\Omega_{j}}(x)dx\\
&=\!\! & \ds\sum_{k=1}^{N}\int_{\Omega}
(\varphi(x)-\varphi(x_{k}))^{2}
\chi_{\Omega_{k}}(x)dx.\nonumber\\
\nonumber
\end{eqnarray}

Next, we  find an estimate for
$$|\varphi(x)-\varphi(x_{k})|^{2}.$$

Consider $\Omega_{k}$ for $k$ fixed, but arbitrary. Choose $z\in\Omega_{k}$ such that $z$ is in the line of the intersection of two planes: the plane which contains the point $x$ and is parallel to $xy$-plane and the plane which contais the point $x_{k}$ and is parallel to the $xz$-plane in three-dimensions.

In other words, if $x$ and $x_{k}$ are such that
$x=(\xi_{1},\xi_{2},\xi_{3})$ and $x_{k}=(\eta_{1},\eta_{2},\eta_{3})$, then $z=(\tau_{1},\eta_{2},\xi_{3})$. Therefore
\begin{eqnarray}
|\varphi(x)-\varphi(x_{k})| & \leq & |\varphi(x)-\varphi(z)|+|\varphi(z)-\varphi(x_{k})|\nonumber\\
& \leq & |\varphi(\xi_{1},\xi_{2},\xi_{3})-\varphi(\tau_{1},\eta_{2},\xi_{3})|+
|\varphi(\tau_{1},\eta_{2},\xi_{3})-\varphi(\eta_{1},\eta_{2},\eta_{3})|.\nonumber\\
\nonumber
\end{eqnarray}
Now me make use of Lemma \ref{difdospontos}, applying $(\ref{terccoordigual})$ for the difference $|\varphi(\xi_{1},\xi_{2},\xi_{3})-\varphi(\tau_{1},\eta_{2},\xi_{3})|$ and $(\ref{segundcoordigual})$ for the difference $|\varphi(\tau_{1},\eta_{2},\xi_{3})-\varphi(\eta_{1},\eta_{2},\eta_{3})|$:

$$\begin{array}{lll}
|\varphi(x)-\varphi(x_{k})|^{2} & \leq &
 (|\varphi(\xi_{1},\xi_{2},\xi_{3})-\varphi(\tau_{1},\eta_{2},\xi_{3})|+
|\varphi(\tau_{1},\eta_{2},\xi_{3})-\varphi(\eta_{1},\eta_{2},\eta_{3})|)^{2}\\
\\
& \leq & 2|\varphi(\xi_{1},\xi_{2},\xi_{3})-\varphi(\tau_{1},\eta_{2},\xi_{3})|^{2}+2|\varphi(\tau_{1},\eta_{2},\xi_{3})-\varphi(\eta_{1},\eta_{2},\eta_{3})|^{2}\\
\\
& \leq &
\ds\frac{4}{h}\left(
4\|\nabla\varphi\|^{2}_{L^{2}(\Omega_{k})}+h^{2}
\left\|\frac{\partial^{2}\varphi}{\partial x\partial y}\right\|^{2}_{L^{2}(\Omega_{k})}\right)+\ds\frac{4}{h}\left(
4\|\nabla\varphi\|^{2}_{L^{2}(\Omega_{k})}+h^{2}
\left\|\frac{\partial^{2}\varphi}{\partial x\partial z}\right\|^{2}_{L^{2}(\Omega_{k})}\right),\\
\end{array}$$
where $h$ is the edge of de cubes $\Omega_{k}$, ie, $h=L/\sqrt[3]{N}$. Then we conclude that
\begin{equation}\label{difdospontosestimada}
|\varphi(x)-\varphi(x_{k})|^{2}\leq
\ds\frac{4}{h}\left(8\|\nabla\varphi\|^{2}_{L^{2}(\Omega_{k})}+h^{2}\left\|\frac{\partial^{2}\varphi}{\partial x\partial y}\right\|^{2}_{L^{2}(\Omega_{k})}+h^{2}
\left\|\frac{\partial^{2}\varphi}{\partial x\partial z}\right\|^{2}_{L^{2}(\Omega_{k})}\right)
\end{equation}
Therefore from $(\ref{somatorianodal1})$ and $(\ref{difdospontosestimada})$, it follows that
\begin{eqnarray}
|\varphi & \!\!\!\! - & \!\!\!\!\ds\sum_{k=1}^{N}\varphi(x_{k})\chi_{\Omega_{k}}|^{2}\leq
 \ds\sum_{k=1}^{N}\int_{\Omega}
(\varphi(x)-\varphi(x_{k}))^{2}
\chi_{Q_{k}}(x)dx\nonumber\\
& \leq & \sum_{k=1}^{N}\ds\int_{\Omega}\ds\frac{4}{h}\left(8\|\nabla\varphi\|^{2}_{L^{2}(\Omega_{k})}+h^{2}\left\|\frac{\partial^{2}\varphi}{\partial x\partial y}\right\|^{2}_{L^{2}(\Omega_{k})}+h^{2}
\left\|\frac{\partial^{2}\varphi}{\partial x\partial z}\right\|^{2}_{L^{2}(\Omega_{k})}\right)\chi_{\Omega_{k}}dx\nonumber\\
& = & \sum_{k=1}^{N}\left(\ds\frac{32}{h}\|\nabla\varphi\|^{2}_{L^{2}(\Omega_{k})}\ds\int_{\Omega}\chi_{\Omega_{k}}dx+4h\left\|\frac{\partial^{2}\varphi}{\partial x\partial y}\right\|^{2}_{L^{2}(\Omega_{k})}\ds\int_{\Omega}\chi_{\Omega_{k}}dx+4h
\left\|\frac{\partial^{2}\varphi}{\partial x\partial z}\right\|^{2}_{L^{2}(\Omega_{k})}\ds\int_{\Omega}\chi_{\Omega_{k}}dx\right)\nonumber\\
\nonumber
\end{eqnarray}
Since $|\Omega_{k}|=h^{3}$ for all $k=1,...,N$, we obtain
$$|\varphi -\sum_{k=1}^{N}\varphi(x_{k})\chi_{Q_{k}}|^{2}\leq 32h^{2}\|\nabla\varphi\|^{2}_{L^{2}(\Omega)}+4h^{4}\left\|\frac{\partial^{2}\varphi}{\partial x\partial y}\right\|^{2}_{L^{2}(\Omega)}+4h^{4}\left\|\frac{\partial^{2}\varphi}{\partial x\partial z}\right\|^{2}_{L^{2}(\Omega)},$$
and thus
\begin{equation}
|\varphi -\sum_{k=1}^{N}\varphi(x_{k})\chi_{Q_{k}}|^{2}
\leq 32h^{2}\|\nabla\varphi\|^{2}_{L^{2}(\Omega)}+8h^{4}\|\varphi\|^{2}_{H^{2}(\Omega)}.
\end{equation}

\Bo

\noindent{\bf Acknowledgements}

E.S.T. is thankful to  the kind hospitality of the Universidade Federal do Rio de Janeiro (UFRJ) and Instituto Nacional de Matem\' {a}tica  Pura  e Aplicada (IMPA) where part of this work was completed.
The work of H.J.N.L. is supported in part by CNPq grant \# 306331 / 2010-1 and FAPERJ grant \# E-26/103.197/2012. The work of  E.S.T.  is supported in part by the NSF grants  DMS-1009950, DMS-1109640 and DMS-1109645, as well as by the  CNPq-CsF grant \# 401615/2012-0, through the program Ci\^encia sem Fronteiras.

\end{document}